\documentclass[a4paper,oneside,11pt]{amsart}

\reversemarginpar

\usepackage[colorlinks,hyperindex,linkcolor=blue,urlcolor=black,pdftitle={One dimensional perturbations of unitaries that are quasiaffine transforms of singular unitaries, and 
multipliers between model spaces},pdfauthor={Maria F.Gamal'}]
{hyperref}

\usepackage{amsmath}
\usepackage{amsfonts}
\usepackage{amssymb}
\usepackage{amsthm}

\usepackage[english]{babel}
\usepackage{enumerate}

\usepackage{graphicx}
\usepackage{color}

\usepackage[abbrev]{amsrefs}

\numberwithin{equation}{section}

\theoremstyle{plain}
\newtheorem{theorem}{Theorem}[section]
\newtheorem{lemma}[theorem]{Lemma}
\newtheorem{corollary}[theorem]{Corollary}
\newtheorem{proposition}[theorem]{Proposition}

\theoremstyle{definition}
\newtheorem{remark}[theorem]{Remark}
\newtheorem{example}[theorem]{Example}

\begin{document}

\title[One dimensional perturbations of  unitaries]{One dimensional perturbations of  unitaries that are quasiaffine transforms of singular unitaries, and 
multipliers between model spaces}

\author{Maria F. Gamal'}
\address{
 St. Petersburg Branch\\ V. A. Steklov Institute 
of Mathematics\\
 Russian Academy of Sciences\\ Fontanka 27, St. Petersburg\\ 
191023, Russia  
}
\email{gamal@pdmi.ras.ru}





\begin{abstract}
 It is shown that, under some natural additional conditions, an operator which intertwines one cyclic  singular unitary operator with 
one dimensional perturbation of another cyclic  singular  unitary operator is the operator of multiplication by a multiplier between model spaces. 
Using this result, it is shown that if $T$ is one dimensional perturbation of a unitary operator, $T$  is a quasiaffine transform of a 
singular unitary operator, and $T$ is power bounded, then $T$ is similar to a unitary operator. Moreover,  
$\sup_{n\geq 0}\|T^{-n}\|\leq(2(\sup_{n\geq 0}\|T^n\|)^2+1)\cdot(\sup_{n\geq 0}\|T^n\|)^5$.
  \end{abstract}

\maketitle

\section{Introduction}

Let $\mathcal H$ and $\mathcal K$ be (complex, separable) Hilbert spaces, and let $\mathcal L(\mathcal H, \mathcal K)$ be the space of (linear, bounded) 
operators acting from $\mathcal H$ to $\mathcal K$. If $\mathcal H =\mathcal K$, we will write $\mathcal L(\mathcal H)$ instead of $\mathcal L(\mathcal H, \mathcal H)$.
Let $T\in\mathcal L(\mathcal H)$, and let $R\in\mathcal L(\mathcal K)$. An operator $X\in\mathcal L(\mathcal H, \mathcal K)$ {\it intertwines} $T$ with $R$, if 
 $XT=RX$. If $X$ is a unitary operator, then $T$ and $R$ 
are called {\it unitarily equivalent}, in notation: $T\cong R$. If $X$ is invertible, that is, $X^{-1}\in\mathcal L(\mathcal K, \mathcal H)$, 
then $T$ and $R$ are called {\it similar}, in notation: $T\approx R$.
If $X$ is a {\it quasiaffinity}, that is, $\ker X=\{0\}$ and $\operatorname{clos}X\mathcal H=\mathcal K$, then
$T$ is called a {\it quasiaffine transform} of $R$, in notation: $T\prec R$. Clearly, $T\prec R$ if and only if $R^\ast\prec T^\ast$.
 If $T\prec R$ and $R\prec T$, then $T$ and $R$ are called 
{\it quasisimilar}, in notation: $T\sim R$. Recall that if $T$ and $R$ are unitary operators and $T\prec R$, 
then $T\cong R$ {\cite[Proposition II.3.4]{sfbk}}. Also, $T$ is similar to a unitary operator if and only if $T$ is invertible and 
$\sup_{n\in\mathbb Z}\|T^n\| < \infty$ \cite{szn47}.

An operator $T\in\mathcal L(\mathcal H)$ is called {\it power bounded}, if 
$\sup_{n\geq 0}\|T^n\| < \infty$.  In \cite{ke}  
a question is raised: {\it  if a power bounded operator is quasisimilar to a singular unitary operator, is it similar to a unitary operator?} 
In some particular cases, an answer is positive \cite{r}, \cite{at},  \cite{g}. In particular, if $U$ is a  unitary operator, 
 $W$ is an invertible operator,  there exists a subsequence $\{n_k\}_k$ of positive integers such that $U^{n_k}\to_kW$ weakly,
 $T$ is a power bounded operator,  and $T\prec U$, 
then $$\sup_{n\geq 0}\|T^{-n}\|\leq \|W^{-1}\|(\sup_{n\geq 0}\|T^n\|)^2.$$ Therefore,  by  \cite{szn47}, 
$T$ is similar to a unitary operator, and by {\cite[Proposition II.3.4]{sfbk}}, $T\approx U$. But the author does not know  any example of a power bounded
 operator $T$ satisfying the above conditions and such that $\sup_{n\geq 0}\|T^{-n}\|$ actually depends on $\|W^{-1}\|$. 

In the present paper, 
it is proved that if $T$ is a power bounded operator, $U$ is a singular unitary operator, $T\prec U$ and $T$ is  one dimensional perturbation 
of a unitary operator, then $$\sup_{n\geq 0}\|T^{-n}\|\leq\bigl(2(\sup_{n\geq 0}\|T^n\|)^2+1\bigr)\cdot(\sup_{n\geq 0}\|T^n\|)^5.$$ Consequently, 
$T\approx U$. If $U$ is an absolutely continuous unitary operator, then there exists a {\it contraction} $T$, that is, $\|T\|\leq 1$, such that $T$ is 
  one dimensional perturbation of a unitary operator, $T\sim U$, and $T\not\approx U$. Also, for every $a>0$ there exists a contraction $T$, 
 such that $T$ is   one dimensional perturbation of a unitary operator, $T\approx U$, and $\|T^{-1}\|>a$. Examples are contractions 
with scalar outer characteristic functions  {\cite[\S VI.3, \S VI.4]{sfbk}}. On the other hand, if $T$ is a 
contraction and $T\prec U$ for a singular unitary operator $U$, then $T\cong U$  
{\cite[Theorems I.3.2, II.6.4, II.2.3]{sfbk}}. 

Clearly, $T$ is one dimensional perturbation of a unitary operator if and only if $T^\ast$ is one dimensional perturbations of a unitary 
operator, and $U$ is a singular unitary operator if and only if $U^\ast$ is a singular unitary operator. Therefore, one can consider 
one dimensional perturbation $T$ of a unitary operator such that $U\prec T$ for a singular unitary operator $U$.  In the present paper, 
it is shown that the question about  similarity of such $T$ and $U$  
 can be reduced to the case where 
the intertwining operator which realizes the relation $U\prec T$ is the operator of multiplication by a {\it  multiplier between model spaces}. 
Multipliers between model spaces are studed in  \cite{c} and  \cite{fhr}. In particular, in  \cite{c} and  \cite{fhr}  isometric multipliers are studed, that is, 
the operator of multiplication by such multiplier is an isometry. In the case considered in the present paper,  multiplication by an isomertic 
 multiplier realizes the relation $T\cong U$. Therefore, isometric multipliers, as well as such multipliers that the range of 
the operator of multiplication by them is not dense   \cite{fhr}, are not considered in the present paper. 
On the other hand, a result of the present paper can be formulated as follows. {\it  If   the operator $X$ of multiplication by a multiplier realizes the relation $U\prec T$ and $T$ is power bounded, then $T\approx U$.} A natural question arises: {\it  is 
the relation $T\approx U$ realized by $X$?}  As is shown in Theorem~\ref{ttt4} (below), every multiplier with dense range  realizes the relation $U\prec T$ for some operator $T$. As is shown in {\cite[Theorem 6.14]{fhr}}, a multiplier  with dense range  can be  noninvertible. 
But  the relationship between power boundedness of $T$ and invertibility of $X$ does not considered in the present paper.


\medskip

{\bf Main Theorem.} {\it  Suppose $T$ is a power bounded operator, $T$ is one dimensional perturbation of a unitary operator, 
$U$ is a singular unitary operator, and $T\prec U$. 
Then  $$\sup_{n\geq 0}\|T^{-n}\|\leq\bigl(2(\sup_{n\geq 0}\|T^n\|)^2+1\bigr)\cdot(\sup_{n\geq 0}\|T^n\|)^5.$$ Consequently,  
$T\approx U$.}

\medskip

The paper is organized as follows. In Section 2, known facts about model spaces and Clark measures are recalled. In Section 3, it is shown 
that operators of multiplication by multipliers between model spaces intertwine cyclic singular unitary operators with one dimensional 
perturbations of cyclic singular unitary operators. In Section 4, auxiliary results using in Section 5 are proved. In Section 5, it is shown 
that if the relation $U\prec T$ is realized by the operator of multiplication by a multiplier between model spaces and $T$ is invertible, 
 then $\sup_{n\geq 0}\|T^{-n}\|\leq(\sup_{n\geq 0}\|T^n\|)^5$. In Section 6, some results closed to results 
of Section 5 are proved. Namely, asymmetric truncated Toeplitz operators as intertwining operators are considered. In Section 7, Main Theorem is proved.

The following notation is used. For a subspace $\mathcal M$ of a Hilbert space, by $P_{\mathcal M}$ the orthogonal projection on  $\mathcal M$ is denoted. 
 For a positive measure $\nu$, by $U_\nu$ the operator of multiplication by the independent variable in $L^2(\nu)$ is denoted.
The symbols $\mathbb T$, $\mathbb D$,  and $\overline{\mathbb D}$ denote  the unit circle,  the open unit disk, and the closed unit disk,  respectively. The symbols $\chi$ and 
 $\text{\bf 1}$ denote the identity function and the unit constant function, respectively, that is, $\chi(\zeta)=\zeta$, $\text{\bf 1}(\zeta)=1$, $\zeta\in\overline{\mathbb D}$. The symbol $m$ denotes the normalized Lebesgue measure on $\mathbb T$. For 
a finite positive Borel measure $\mu$ on $\mathbb T$ set 
$$\widehat\mu(n)=\int_{\mathbb T}\zeta^{-n}\text{\rm d}\mu(\zeta),  \ \ \ n\in\mathbb Z.$$

\section{ Model spaces and Clark measures}

In this section, we recall facts about model spaces and Clark measures (or Aleksandrov--Clark measures, in other terminology) of inner functions, which will be needed in the sequel, see  \cite{clark} and \cite{polt} for the references, see also {\cite[Sec. 2.7]{gr}}, {\cite[Sec. 8]{grm}} and references therein. 

Denote by $m$ the normalized Lebesgue measure on the unit circle $\mathbb T$ and by $H^2$  the Hardy space in the unit disk $\mathbb D$.
Indeed, functions from $H^2$ have nontangential boundary values on $\mathbb T$ a.e. with respect to $m$. Set $H^2_-=L^2(\mathbb T,m)\ominus H^2$.
Set $\chi(z)=z$, $\text{\bf 1}(z)=1$, $z\in\overline{\mathbb D}$.

Let $\theta$ be an inner function. Put $\mathcal K_\theta=H^2\ominus\theta H^2$. The space $\mathcal K_\theta$ is called the {\it model space} 
corresponding to $\theta$. It is well known and easy to see that $f\in\mathcal K_\theta$ if and only if 
$\theta\overline\chi\overline f\in\mathcal K_\theta$. For $\lambda\in\mathbb D$ 
put 
\begin{equation}
\label{y1}
k_{\theta,\lambda}(z)=\frac{1-\overline{\theta(\lambda)}\theta(z)}{1-\overline\lambda z} \ \ \text{ and } \ \ 
k_{\ast\theta,\lambda}(z)=\frac{\theta(z)-\theta(\lambda)}{z-\lambda}, \ \ \ z\in\mathbb D. 
\end{equation}
Then $k_{\theta,\lambda}$, $k_{\ast\theta,\lambda}\in\mathcal K_\theta$, $k_{\ast\theta,\lambda}=\theta\overline\chi\overline{k_{\theta,\lambda}}$,
and $f(\lambda)=(f,k_{\theta,\lambda})$ for any $\lambda\in\mathbb D$ and $f\in\mathcal K_\theta$.

The space $\mathcal K_\theta$ is coinvariant subspace of $S\in\mathcal L(H^2)$, 
the operator of multiplication by $\chi$ in $H^2$. Put $S_\theta=P_{\mathcal K_\theta}S|_{\mathcal K_\theta}$. It is easy to see that 
 $\text{\bf 1}$,  $\overline\chi\theta\in\mathcal K_\theta$ if and only if $\theta(0)=0$. 

Let $\theta(0)=0$. Then 
\begin{equation}
\label{y2}
S_\theta f =\chi f-(f,\overline\chi\theta)\theta, \ \ \ f\in\mathcal K_\theta. 
\end{equation}
For $c\in\mathbb T$ put 
\begin{equation}
\label{y3}
U_{(\theta)c}=S_\theta+c(\cdot,\overline\chi\theta)\text{\bf 1}.\end{equation}
It is easy to see that $U_{(\theta)c}$ is unitary for every $c\in\mathbb T$. 

Let $\mu$ be  a positive Borel singular measure on $\mathbb T$, and let $\mu(\mathbb T)=1$. Then the function $\theta$ defined by the formula
\begin{equation}
\label{yy0}
\frac{1}{1-\theta(z)}=\int_{\mathbb T}\frac{1}{1-z\overline\zeta}\text{\rm d}\mu(\zeta),  \ \ z\in\mathbb D,\end{equation}
is inner, and $\theta(0)=0$. Furthermore, $\theta$ has nontangential boundary values equal to $1$ on $\mathbb T$ a.e. with respect to $\mu$. 
Conversely, for every inner function $\theta$ such that $\theta(0)=0$ there exists 
a positive Borel singular measure $\mu:=\sigma_1$ on $\mathbb T$ such that $\mu(\mathbb T)=1$ and (\ref{yy0}) is fulfilled. Let $c\in\mathbb T$. 
Applying (\ref{yy0}) to $\overline c\theta$ one obtains a  positive Borel singular measure $\sigma_c$ on $\mathbb T$ such that $\sigma_c(\mathbb T)=1$ and 
\begin{equation}
\label{yyc}\frac{1}{1-\overline c\theta(z)}=\int_{\mathbb T}\frac{1}{1-z\overline\zeta}\text{\rm d}\sigma_c(\zeta),  \ \ z\in\mathbb D.
\end{equation}
Since $\theta$ has nontangential boundary values equal to $c$ on $\mathbb T$ a.e. with respect to $\sigma_c$, we obtain that the measures 
$\sigma_c$ are pairwise singular. Every function $f\in\mathcal K_\theta$ has nontangential boundary values $f(\zeta)$ for 
a.e. $\zeta\in\mathbb T$ with respect to $\sigma_c$, for all $c\in\mathbb T$ \cite{polt}. The operator
\begin{equation}
\label{y4}\begin{aligned}&
J_{\theta,c}\in\mathcal L(\mathcal K_\theta,L^2(\mathbb T,\sigma_c)), \\ &  (J_{\theta,c}f)(\zeta)=f(\zeta) \ \ \text{ for a.e. }\ \zeta\in\mathbb T 
\text{ with respect to }  \sigma_c, \ \   f\in\mathcal K_\theta, \end{aligned}\end{equation}
is unitary, and 
\begin{equation}
\label{yy1}
 (J_{\theta,c}^{-1}\gamma)(z)=(1-\overline c\theta(z))\int_{\mathbb T}\frac{\gamma(\zeta)}{1-z\overline\zeta}\text{\rm d}\sigma_c(\zeta),  
\ \ z\in\mathbb D,\ \ \gamma\in L^2(\sigma_c)\end{equation}
\cite{clark}. Furthermore, by (\ref{y2}) and (\ref{y3}), 
\begin{equation}
\label{y5}J_{\theta,c}U_{(\theta)c}=U_{\sigma_c} J_{\theta,c}. \end{equation}

\begin{lemma}
\label{lll21} Suppose $\theta$ is an inner function, $\theta(0)=0$, $u\in\mathcal K_\theta$, and  
 \begin{equation}
T=S_\theta+(\cdot, \overline\chi\theta)u.\end{equation}
If $T$ is an isometry, then $u=c$ a.e. on $\mathbb T$ for some $c\in\mathbb T$. Conversely, if  $u=c$ a.e. on $\mathbb T$  for some $c\in\mathbb T$, 
 then $T$ is unitary.
\end{lemma}

\begin{proof} Indeed,  $T$ is an isometry if and only if $T^\ast T=I$.  It is easy to compute that 
\begin{equation*}
T^\ast T=I-(\cdot,\overline\chi\theta)\overline\chi\theta+(\cdot,\overline\chi\theta)S_\theta^\ast u +
(\cdot,S_\theta^\ast u)\overline\chi\theta+\|u\|^2(\cdot,\overline\chi\theta)\overline\chi\theta. \end{equation*}
If $T^\ast T=I$, then \begin{equation}\label{3.9}\overline\chi\theta=S_\theta^\ast u +
(\overline\chi\theta,S_\theta^\ast u)\overline\chi\theta+\|u\|^2\overline\chi\theta. \end{equation}
Therefore, $S_\theta^\ast u = a\overline\chi\theta$ for some $a\in\mathbb C$. But $\overline\chi\theta\perp S_\theta^\ast\mathcal K_\omega$, 
and we conlcude that $a=0$. Thus, $u=c$ a.e. on $\mathbb T$  for some $c\in\mathbb C$. Clearly, $\|u\|^2=|c|^2$. 
We infer from (\ref{3.9}) that $\|u\|^2=1$. 

Now suppose that $u=c$ a.e. on $\mathbb T$  for some $c\in\mathbb T$. Then,  by (\ref{y3}), $T=U_{(\theta)c}$. Thus,  $T$ is unitary. 
\end{proof}

\smallskip

The following lemma will be used in the proof of Proposition  \ref{ppp1}.

\begin{lemma}
\label{lll7} Suppose $\mu$ is a positive Borel singular measure on $\mathbb T$, $\mu(\mathbb T)=1$, 
$J_{\theta,1}^{-1}\in\mathcal L(L^2(\mu), \mathcal K_\theta)$ acts by the formula (\ref{yy1}),  $\chi(\zeta)=\zeta$, $\zeta\in\overline{\mathbb D}$, 
 $\gamma\in L^2(\mu)$, and $a=\int_{\mathbb T} \gamma \text{\rm d}\mu$.
Set $u=J_{\theta,1}^{-1}(\overline\chi\overline\gamma)$. Then $$J_{\theta,1}^{-1}\gamma = \theta \overline\chi\overline u \ \ \text{ and }\ \  
J_{\theta,1}^{-1}\overline\gamma = \chi u +\overline a(1-\theta).$$
\end{lemma}
\begin{proof} For $\lambda\in\mathbb D$ set $\gamma_\lambda(\zeta)=\frac{1-\overline{\theta(\lambda)}}{1-\overline\lambda\zeta}$, $\zeta\in\mathbb T$. 
We have $J_{\theta,1}k_{\theta,\lambda}=\gamma_\lambda$ and  $J_{\theta,1}k_{\ast\theta,\lambda}=\overline\chi\overline\gamma_\lambda$. Thus, for 
$u=J_{\theta,1}^{-1}\overline\chi\overline\gamma_\lambda$ we have $\theta\overline\chi\overline u = J_{\theta,1}^{-1}\gamma_\lambda$. Since 
$\{\gamma_\lambda\}_{\lambda\in\mathbb D}$ is total in $L^2(\mu)$, we have that $J_{\theta,1}^{-1}\gamma = \theta \overline\chi\overline u$, where 
$u=J_{\theta,1}^{-1}(\overline\chi\overline\gamma)$, for every $\gamma\in L^2(\mu)$. 

Recall that $U_\mu\in\mathcal L(L^2(\mu))$ is the operator of multiplication by $\chi$. 
We have $$J_{\theta,1}^{-1}\overline\gamma = J_{\theta,1}^{-1}U_\mu (\overline\chi\overline\gamma) = J_{\theta,1}^{-1}U_\mu J_{\theta,1}J_{\theta,1}^{-1}(\overline\chi\overline\gamma) = 
J_{\theta,1}^{-1}U_\mu J_{\theta,1}u.$$
Since $J_{\theta,1}$ is unitary, 
 $$(u,\overline\chi\theta)_{\mathcal K_\theta}=(J_{\theta,1}^{-1}(\overline\chi\overline\gamma),J_{\theta,1}^{-1}\overline\chi)_{\mathcal K_\theta}=
(\overline\chi\overline\gamma,\overline\chi)_{L^2(\mu)}=\overline a.$$
By (\ref{y5}), (\ref{y2}), and (\ref{y3}), 
$$(J_{\theta,1}^{-1}U_\mu J_{\theta,1}f)(z)=zf(z)+(f,\overline\chi\theta)_{\mathcal K_\theta}(1-\theta(z)), \ \ \  z\in\mathbb D, 
 \ \ \ f\in \mathcal K_\theta.$$
Therefore, \begin{equation*}\begin{aligned}(J_{\theta,1}^{-1}\overline\gamma)(z)=(J_{\theta,1}^{-1}U_\mu J_{\theta,1}u)(z) & =zu(z)+(u,\overline\chi\theta)_{\mathcal K_\theta}(1-\theta(z))\\ & =
zu(z)+\overline a(1-\theta(z)), \ \  z\in\mathbb D. \end{aligned}\end{equation*}
\end{proof}

\section{ Multipliers between model spaces as intertwining operators}

In this section, it is proved that, under some natural additional conditions, an operator which intertwines a cyclic singular unitary operator with one dimensional perturbation of another cyclic singular unitary operator, is the operator of multiplication by a multiplier between model spaces. Also, the converse is true.

Suppose $\theta$, $\omega\in H^\infty$ are inner functions, and $g\in H^2$. If  $gf\in \mathcal  K_\omega$ for every $f\in\mathcal K_\theta$, 
then $g$ is called a {\it multiplier} between $\mathcal K_\theta$ and $\mathcal K_\omega$. Clearly, the mapping 
\begin{equation}\label{dddmult} \mathcal K_\theta\to \mathcal K_\omega, \ \ \ f\mapsto gf,  \ \ \ f\in\mathcal K_\theta,\end{equation} 
is bounded by the closed graph theorem.

The following facts can be found in \cite{c}, \cite{fhr}, or can be easy checked straightforward. 

\begin{lemma}\label{lllmult} \cite{c}, \cite{fhr}.  Suppose $\theta$, $\omega\in H^\infty$ are inner functions, and $g\in H^2$. Set 
$\chi(z)=z$, $z\in\mathbb T$.

(i)  The function $g$ is  a  multiplier between $\mathcal K_\theta$ and $\mathcal K_\omega$ if and only if there exists a function 
$g_1\in H^2$ such that $g=\omega\overline\theta\overline g_1$ a.e. on $\mathbb T$ and $gf\in H^2$ for every $f\in\mathcal K_\theta$. 

(ii) If $\operatorname{clos}g\mathcal K_\theta=\mathcal K_\omega$, then $g$ is outer, and  $g=c\omega\overline\theta\overline g$ a.e. on $\mathbb T$ for some $c\in\mathbb T$.  Furthermore,   $g\mathcal K_\theta=\mathcal K_\omega$ if 
and only if  $f/g\in H^2$ for every $f\in\mathcal K_\omega$.  
\end{lemma}

\begin{remark} Crofoot transform allows to consider Frostman shift of inner functions $\theta$ and $\omega$ instead of 
 $\theta$ and $\omega$ themselves when multipliers are studed (see   \cite{c}, \cite{fhr}, {\cite[Sec. 3.3]{grm}}). Therefore, we can  assume that $\theta(0)=0$ and $\omega(0)$ =0. \end{remark} 

\begin{example}\cite{c}, \cite{fhr} (i) There exist isometric multipliers $g$, i.e.,  such  that   the operator of multiplication by $g$ defined by (\ref{dddmult}) is a unitary operator \cite{c}.

(ii) There exist nonisometric multipliers $g$ such  that $g\mathcal K_\theta=\mathcal K_\omega$,
 in other words,  the operator of multiplication by $g$ defined by (\ref{dddmult}) is nonunitary invertible operator \cite{c}.

(iii) There exist multipliers $g$ such  that $g\mathcal K_\theta\neq\mathcal K_\omega$ and  the operator of multiplication by $g$ defined by (\ref{dddmult}) is left-invertible. Namely, let $\theta$ and $\vartheta$ be two nonconstant  inner functions such that their boundary spectrums do not intersect (see \cite{fhr} for definition). By {\cite[Theorem 3.5]{fhr}}, the set of multipliers between $\mathcal K_\theta$ and 
$\mathcal K_{\vartheta\theta}$ is $\mathcal K_{\chi\vartheta}$.  Take $g\in\mathcal K_{\chi\vartheta}$ such that $1/g\in H^\infty$ 
(for example, $g$  can be  defined by (\ref{y1}) applied to $\chi\vartheta$). Then the operator of multiplication by $g$ is left-invertible. 
If this operator is onto $\mathcal K_{\vartheta\theta}$, then, by \cite{c}, $\dim\mathcal K_{\chi\vartheta}=1$, a contradiction with 
the assumption that $\vartheta$ is nonconstant.

(iv) There exist multipliers $g$ such  that $\operatorname{clos}g\mathcal K_\theta=\mathcal K_\omega$ and $g\mathcal K_\theta\neq\mathcal K_\omega$  {\cite[Theorem 6.14]{fhr}} (the first equality follows from condition (ii) in the proof of
 {\cite[Theorem 6.14]{fhr}}).
\end{example}

  The following lemma can be easy checked straightforward, therefore, its proof is omitted. 

  \begin{lemma}
\label{lll10} Suppose $R\in \mathcal L(\mathcal H)$, $u$, $v\in \mathcal H$, $T=R+(\cdot, v)u$, and $\lambda\in\mathbb C$. 
If $(R-\lambda)^{-1}$ exists and $1+((R-\lambda)^{-1}v, u)\neq 0$, then 
$$(T-\lambda)^{-1}x=(R-\lambda)^{-1}x-\frac{((R-\lambda)^{-1}x, v)}{1+((R-\lambda)^{-1}u, v)}(R-\lambda)^{-1}u, \ \ \ x\in\mathcal H.$$
\end{lemma}

The following theorem is the main result of this section.

\begin{theorem}
\label{ttt3} Suppose $\theta$, $\omega\in H^\infty$ are inner functions, $\theta(0)=0$, $\omega(0)=0$, $\text{\bf 1}(z)=1$, 
$\chi(z)=z$, $z\in\mathbb T$, 
\begin{equation*}\begin{gathered}u\in\mathcal K_\omega,  \ \ \ \ T=S_\omega+(\cdot, \overline\chi\omega)u, \\ X\in \mathcal L(\mathcal K_\theta,\mathcal K_\omega), \ \  
X^\ast\overline\chi\omega = \overline\chi\theta, \ \text{ and } \ XU_{(\theta)1}=TX.\end{gathered}\end{equation*} 
Set $$g=\frac{u-\omega}{1-\theta}.$$ Then $g\in\mathcal K_\omega$ and $Xf=gf$ for every $f\in\mathcal K_\theta$.
\end{theorem}
\begin{proof} Recall that  $k_{\theta,\lambda}$ is defined in (\ref{y1}). We have  
$$((S_\theta-\lambda)^{-1})^\ast\overline\chi\theta=-\frac{1}{\overline{\theta(\lambda)}}k_{\theta,\lambda}$$
for every $\lambda\in\mathbb D$ such that $\theta(\lambda)\neq 0$. 
Therefore, 
\begin{equation}
\label{3.1}((S_\theta-\lambda)^{-1}f,\overline\chi\theta)=-\frac{f(\lambda)}{\theta(\lambda)}
\end{equation}
for every $f\in\mathcal K_\theta$ and every $\lambda\in\mathbb D$ such that $\theta(\lambda)\neq 0$.
 
Let $f\in\mathcal K_\theta$, and let $\lambda\in\mathbb D$ be such that $\theta(\lambda)\neq 0$, $\omega(\lambda)\neq 0$, and 
$u(\lambda)\neq \omega(\lambda)$.
By Lemma \ref{lll10} applied to $T$ and (\ref{3.1}) applied to $\omega$, 
\begin{equation}
\label{3.2}
\begin{aligned}((T-&\lambda)^{-1}Xf, \overline\chi\omega)\\&=((S_\omega-\lambda)^{-1}Xf,\overline\chi\omega)-
\frac{((S_\omega-\lambda)^{-1}Xf, \overline\chi\omega)}{1+((S_\omega-\lambda)^{-1}u, \overline\chi\omega)}
((S_\omega-\lambda)^{-1}u,\overline\chi\omega)\\&
=
-\frac{(Xf)(\lambda)}{\omega(\lambda)}-\frac{(Xf)(\lambda)}{\omega(\lambda)}\frac{1}{1-\frac{u(\lambda)}{\omega(\lambda)}}
\frac{u(\lambda)}{\omega(\lambda)}=\frac{(Xf)(\lambda)}{u(\lambda)-\omega(\lambda)}.
\end{aligned}
\end{equation}
By Lemma \ref{lll10} applied to $U_{(\theta)1}$ and (\ref{3.1}) applied to $\theta$,  
\begin{equation}
\label{3.3}
\begin{aligned}&\Bigl(X\Bigl((S_\theta-\lambda)^{-1}f-\frac{((S_\theta-\lambda)^{-1}f, \overline\chi\theta)}{1+((S_\theta-\lambda)^{-1}\text{\bf 1}, 
\overline\chi\theta)}
(S_\theta-\lambda)^{-1}\text{\bf 1}\Bigr),\overline\chi\omega\Bigr)\\&
=\Bigl
((S_\theta-\lambda)^{-1}f-\frac{((S_\theta-\lambda)^{-1}f, \overline\chi\theta)}{1+((S_\theta-\lambda)^{-1}\text{\bf 1}, \overline\chi\theta)}
(S_\theta-\lambda)^{-1}\text{\bf 1},\overline\chi\theta\Bigr)=\frac{f(\lambda)}{1-\theta(\lambda)}.
\end{aligned}
\end{equation}
Since $$X(U_{(\theta)1}-\lambda)^{-1}=(T-\lambda)^{-1}X,$$ we obtain from 
(\ref{3.2}) and (\ref{3.3}) that $$(Xf)(\lambda)=\frac{u(\lambda)-\omega(\lambda)}{1-\theta(\lambda)}f(\lambda)$$
for every $f\in\mathcal K_\theta$ and every $\lambda\in\mathbb D$ such that  $\theta(\lambda)\neq 0$, $\omega(\lambda)\neq 0$, and 
$u(\lambda)\neq \omega(\lambda)$. Therefore, $Xf=gf$. Since $\text{\bf 1}\in\mathcal K_\theta$ and 
$X\text{\bf 1}=g$, we conclude that $g\in\mathcal K_\omega$. 
\end{proof}

\begin{theorem}\label{ttt4}  Suppose $\theta$, $\omega\in H^\infty$ are inner functions, $\theta(0)=0$, $\omega(0)=0$, 
 $\text{\bf 1}(z)=1$, 
$\chi(z)=z$, $z\in\mathbb T$,  and $g$ is a multiplier between $\mathcal K_\theta$ and $\mathcal K_\omega$.  Let the function $g_1$ be defined in Lemma~\ref{lllmult}(i), and let the operator $X\in \mathcal L(\mathcal K_\theta,\mathcal K_\omega)$ be defined by (\ref{dddmult}). 
 Suppose  $g_1(0)\neq 0$. Let $c\in\mathbb T$. Put 
\begin{equation}
\label{3.4}
 u=\omega+\frac{1}{\overline {g_1(0)}}(c-\theta)g.\end{equation}
 Then $u\in \mathcal K_\omega$. Put 
\begin{equation}
\label{3.5}T=S_\omega+(\cdot, \overline\chi\omega)u.\end{equation}
 Then \begin{equation}
\label{3.6}XU_{(\theta)c}=TX.\end{equation}
\end{theorem} 
\begin{proof}
 Clearly, $u\in H^2$, 
$$u=\omega\cdot\Bigl(\overline{1+\frac{(\overline c\theta-1)g_1}{g_1(0)}}\Bigr) \text{ a.e. on $\mathbb T$, \ \ \  and } \ \ \ 
\Bigl(1+\frac{(\overline c\theta-1)g_1}{g_1(0)}\Bigr)(0)=0,$$
because $\theta(0)=0$. Thus, $u\in\mathcal K_\omega$.  Furthermore,  
\begin{equation}
\label{3.7new}  X^\ast(\overline\chi\omega)=g_1(0)\overline\chi\theta.   \end{equation}  
It follows from  (\ref{y2}) and (\ref{3.7new}) that 
 \begin{equation}
\label{3.7}S_\omega X-X S_\theta=(\cdot, \overline\chi\theta)\cdot(\theta g - \overline {g_1(0)}\omega).\end{equation}
Now (\ref{3.6}) follows from (\ref{y3}),  (\ref{3.4}), (\ref{3.7new})  and (\ref{3.7}).
\end{proof}

\begin{example}
\label{eee3.7}
\cite{c}. Suppose $\theta$ is an inner function, $\theta(0)=0$, and $0\neq\lambda\in\mathbb D$. Recall that $k_{\theta, \lambda}$ and 
$k_{\ast\theta,\lambda}$ are defined in (\ref{y1}). Put $$ \omega=\chi  k_{\ast\theta,\lambda} / k_{\theta,\lambda} \ \ \ \text{ and } 
\ \ \ g=1/k_{\theta, \lambda}.$$
Clearly, $g$ is outer, and $g=\omega\overline\theta\overline g$ a.e. on $\mathbb T$.  Let $c\in\mathbb T$. Define $u$ by (\ref{3.4}) and $T$ by (\ref{3.5}).
Then  $u=\omega+(c-\theta) g$ and $T\approx U_{(\theta)c}$. If 
\begin{equation}\label{3.10}\operatorname{card}\{w\in\mathbb D\ :\ \theta(w)=\theta(\lambda)\}>2,\end{equation} then $T$ is not unitary. 

Indeed,  $g$ is  a  multiplier between $\mathcal K_\theta$ and $\mathcal K_\omega$ by Lemma~\ref{lllmult}(i). Define $X$ by (\ref{dddmult}). By Theorem \ref{ttt4}, $XU_{(\theta)c}=TX$. By Lemma \ref{lllmult}(i), $X$ is invertible. Thus, 
$X$ realizes the relation $T\approx U_{(\theta)c}$.

 Suppose that $T$ is unitary. By Lemma \ref{lll21},  $u=c_1\text{\bf 1}$  
for some $c_1\in\mathbb T$. Thus, \begin{equation}\label{3.11}c_1=\omega+(c-\theta) g.\end{equation} Take $w\in\mathbb D$ 
such that $w\neq\lambda$ and  $\theta(w)=\theta(\lambda)$. By the definition of $\omega$, $\omega(w)=0$. We infer from (\ref{3.11}) that 
$$w=\frac{1}{\overline\lambda}\Bigl(1-c_1\frac{1-|\theta(\lambda)|^2}{c-\theta(\lambda)}\Bigr).$$
This contradicts with (\ref{3.10}).
\end{example}

In the following example an operator is constructed, which is two dimensional perturbation of one unitary operator, is similar to this unitary operator, and is one dimensional perturbation of another  unitary operator. The idea from  {\cite[Lemma 3.1]{k}} is used. In  {\cite[Lemma 3.1]{k}}, unitary operators are considered, therefore, 
 $|\varphi|=1$ a.e. with respect to $\nu$ is supposed there.  In the present paper, nonunitary operators are considered. The obtained multiplier is from {\cite[Corollary 7.6]{fhr}}. 

\begin{example}
\label{eeekfhr} {\cite[Lemma 3.1]{k}} {\cite[Corollary 7.6]{fhr}}. Suppose  $\theta$ is an inner function, $\theta(0)=0$, and 
$1\neq c\in\mathbb T$. Recall that $U_{(\theta)1}$ is defined by (\ref{y3}), and $\sigma_c$ is defined by (\ref{yyc}). Set  $\nu=\sigma_c$. Let $\varphi\in L^\infty(\nu)$ be such that  $1/\varphi\in L^2(\nu)$. Without loss of generality, 
we can assume that $ \int_{\mathbb T}\frac{1}{|\varphi|^2}{\text{\rm d}}\nu =1$. Put
$$ Y=\varphi(U_\nu) \ \ \ \text {and } \ \ \  T=U_\nu+(\overline c-1)(\cdot,\overline\chi/\overline\varphi)\varphi.$$
Recall that $J_{\theta,c}$ is defined by (\ref{y4}). Taking into account that 
$$ U_{(\theta)1} - U_{(\theta)c} =  (1-c)(\cdot,\overline\chi\theta)\text{\bf 1}, $$
it is easy to see that 
$$ YJ_{\theta,c}U_{(\theta)1}J_{\theta,c}^{-1}=TY \ \ \ \text {and } \ \ \  
 \dim \operatorname{ran} (T-J_{\theta,c}U_{(\theta)1}J_{\theta,c}^{-1})\leq 2 $$
(see (\ref{y2}), (\ref{y3}), and (\ref{y5})).
Put ${\text{\rm d}}\nu_1 =\frac{1}{|\varphi|^2}{\text{\rm d}}\nu$. Then $\nu_1(\mathbb T)=1$.
Define $Z\in\mathcal L(L^2(\nu), L^2(\nu_1))$  by the formula 
$$Zf=\overline\varphi f, \ \ \  f\in L^2(\nu).$$ It is easy to see that $Z$ is a unitary operator. 
  Define $\omega$ by $\nu_1$ as in  (\ref{yy0}). 
 Put $T_1=J_{\omega,1}^{-1}ZTZ^{-1}J_{\omega,1}$.  Then $T\cong T_1$, and it is easy to see that 
\begin{equation*}T_1=S_\omega+(\cdot,\overline\chi\omega)\bigl((\overline c-1)J_{\omega,1}^{-1}|\varphi|^2+1\bigr).
\end{equation*}
Set $X=\overline c J_{\omega,1}^{-1}ZYJ_{\theta,c}$. It is easy to see that 
$$XU_{(\theta)1}=T_1X  \ \ \ \text {and } \ \ \  X^\ast \overline\chi\omega=\overline\chi\theta.$$
Thus, $T_1$ and $X$ satisfy to the conditions of Theorem~\ref{ttt3} with 
\begin{equation}\label{eeekfhr1} u=(\overline c-1)J_{\omega,1}^{-1}|\varphi|^2+1.\end{equation} 
By Theorem~\ref{ttt3}, $g=\frac{u-\omega}{1-\theta}$ is a multiplier between $\mathcal K_\theta$ and $\mathcal K_\omega$.
Applying (\ref{yy1}) and  (\ref{yyc}),  it is easy to see that 
\begin{equation}\label{eeekfhr2}J_{\omega,1}^{-1}|\varphi|^2=\frac{1-\omega}{1-\overline c\theta}.\end{equation}
Therefore, $$g=\overline c  \frac{1-\omega}{1-\overline c\theta}.$$
Applying Lemma~\ref{lll21} to $T_1$ and taking into account (\ref{eeekfhr1}) and (\ref{eeekfhr2}), it is easy to see that $T_1$ is unitary 
if and only if $\omega=\overline c \theta$. The latest equality means that $\nu_1=\sigma_c=\nu$. Thus, $T_1$ is unitary 
if and only if $|\varphi|=1$ a.e. with respect to $\nu$. Therefore, if $1/\varphi\in L^\infty(\nu)$, and $|\varphi|$ is a nonconstant function,    then $Y$ is invertible, $T_1$ is nonunitary, and  $T_1\approx U_{(\theta)1}$. 
\end{example}

\section{ Auxiliary results}

In this section,  auxiliary results using in the next section are proved.

\smallskip

The following lemma is a particular case of {\cite[Theorem 2.5 and Lemma 2.4]{g}}, see also \cite{at}, \cite{r}. Note that a misprint is in 
{\cite[Theorem 2.5]{g}}: $x$ must be deleted twice.

\begin{lemma}
\label{lll8}  Suppose $\{\lambda_j\}_j\subset \mathbb T$, $\{\mathcal H_j\}_j$ is 
a family of Hilbert spaces,  $$\mathcal H=\oplus_j \mathcal H_j, \ \ \ U=\oplus_j \lambda_jI_{\mathcal H_j},$$
  $R\in\mathcal L(\mathcal K)$ is a  power bounded operator, and $R\prec U$. Then there exists an invertible operator 
$Y\in\mathcal L(\mathcal K,\mathcal H)$ such that $$YR=UY, \ \ \ \|Y\|\leq \sup_{n\geq 0}\|R^n\|, \ \ \text{ and } \ \ 
\|Y^{-1}\|\leq (\sup_{n\geq 0}\|R^n\|)^2.$$
\end{lemma}
\begin{proof}  By {\cite[Theorem 2.5 and Lemma 2.4]{g}}, $R\approx U$ and $$\sup_{n\geq 0}\|R^{-n}\|\leq (\sup_{n\geq 0}\|R^n\|)^2.$$
 Let $Y$ be the canonical interwining mapping of $R$ and its unitary asymptote $U$ (see \cite{ke}). We have 
$$ \liminf_{n\to +\infty}\|R^nx\|\leq\|Yx\|\leq\limsup_{n\to +\infty}\|R^nx\| \ \ \ \text{ for every }\ x \in \mathcal K.$$
The needed estimate on $\|Y\|$  follows from the right inequality. 

Let $ x \in \mathcal K$. Set $y=Yx$. We have $$\|Y^{-1}y\| = \|x\| = \|R^{-n}R^nx\|\leq (\sup_{n\geq 0}\|R^{-n}\|)\cdot\|R^nx\|.$$ 
Therefore, \begin{equation*}\begin{aligned}\|Y^{-1}y\|&\leq  (\sup_{n\geq 0}\|R^{-n}\|)\cdot \liminf_{n\to +\infty}\|R^nx\| \\  &\leq(\sup_{n\geq 0}\|R^{-n}\|)\cdot \|Yx\| = 
(\sup_{n\geq 0}\|R^{-n}\|)\cdot \|y\|.\end{aligned}\end{equation*} Thus, $$\|Y^{-1}\|\leq\sup_{n\geq 0}\|R^{-n}\|\leq  (\sup_{n\geq 0}\|R^n\|)^2.$$
\end{proof}

\begin{lemma}
\label{lll3} Suppose $\{\lambda_j\}_j$, $\{\xi_j\}_j\subset \mathbb T$, $\{\mathcal H_j\}_j$ is 
a family of Hilbert spaces,  $\mathcal H=\oplus_j \mathcal H_j$, $U$, $W\in\mathcal L(\mathcal H)$, 
$$U=\oplus_j \lambda_jI_{\mathcal H_j}, \ \ \ W=\oplus_j \xi_jI_{\mathcal H_j},$$  
 $\{n_k\}_k$ is a  subsequence of positive integers,  and $$U^{n_k}x\to_k Wx \ \ \text{ for every } \ x\in\mathcal H.$$ 
Furthermore, suppose  $T\in\mathcal L(\mathcal K)$, $M=\sup_{n\geq 0}\|T^n\|$, 
$X\in\mathcal L(\mathcal H,\mathcal K)$ is a quasiaffinity, and $$XU=TX.$$ Then there exist $R\in\mathcal L(\mathcal K)$ and 
 an invertible operator $Y\in\mathcal L(\mathcal H,\mathcal K)$ such that 
$$T^{n_k}x\to_k Rx \ \ \text{ for every } \ x\in\mathcal H, \ \ \ XW=RX,$$   $\|Y\|\leq M$, $\|Y^{-1}\|\leq M^2$, 
and $R=YWY^{-1}$. Consequently, $R$ is invertible, $\|R^{-1}x\|\leq M^3\|Rx\|$, and $\|(R^\ast)^{-1}x\|\leq M^3\|R^\ast x\|$ 
for every $x\in\mathcal K$.
\end{lemma}
\begin{proof} Set $\mathcal M_j=\operatorname{clos}X\mathcal H_j$ for every $j$, then $\mathcal K=\vee_j\mathcal M_j$ and 
$T|_{\mathcal M_j}=\lambda_jI_{\mathcal M_j}$. Note that $\lambda_j^{n_k}\to_k\xi_j$ for every $j$. 
Let $\{x_j\}_j\subset\mathcal K$ be such that $x_j\in\mathcal M_j$ for every $j$, and let the cardinality of $\{j\ :\ x_j\neq 0\}$ be finite. 
Put $x=\sum_jx_j$ and $Rx=\sum_j\xi_jx_j$. Since $T^{n_k}x\to_kRx$, we have $\|Rx\|\leq M\|x\|$. Since the set of $x$ constructed as above 
is dense in $\mathcal K$, the mapping $R$ can be extended on the  whole $\mathcal K$, that is, $R\in\mathcal L(\mathcal K)$. 
Since $T^{ln_k}x\to_kR^lx$ for every $l\in\mathbb N$ and every $x$ constructed as above, we obtain that $\|R^lx\|\leq M\|x\|$  for such $x$. 
Since the set of such $x$ is dense in $\mathcal K$, we conclude that $\sup_{n\geq 0}\|R^n\|\leq M$. 
It  easy follows from the definition of $R$ that $RX=XW$. Therefore, $X^\ast R^\ast = W^{-1}X^\ast$. Applying Lemma \ref{lll8} to 
$R^\ast$, 
we obtain an invertible operator $Y_\ast$ such that $Y_\ast R^\ast = W^{-1}Y_\ast$, $\|Y_\ast\|\leq M$, and $\|Y_\ast^{-1}\|\leq M^2$. 
Set $Y=Y_\ast^\ast$. Then  $Y$ satisfies the conclusion of the lemma. Finally, from the equalities 
\begin{equation*}
\begin{aligned} &Y^{-1}=W^{-1}Y^{-1}R \\  \text{ and } \ \  & R^{-1}x=YW^{-1}Y^{-1}x = YW^{-1}W^{-1}Y^{-1}Rx, \ \ \ x\in \mathcal K,\end{aligned}\end{equation*} we obtain that 
$$\|R^{-1}x\|\leq\|Y\|\|W^{-1}\|\|W^{-1}\|\|Y^{-1}\|\|Rx\|\leq M^3\|Rx\|, \ \ x\in \mathcal K.$$ The estimate for $R^\ast$ is proved similarly.
\end{proof}
	
\begin{lemma}
\label{lll1} Suppose $U\in\mathcal L(\mathcal H)$ is a unitary operator, and  $\{\lambda_j\}_j\subset \mathbb T$ are all eigenvalues of $U$. 
Then $\mathcal H=\mathcal H_{\text{\rm at}}\oplus\mathcal H_{\text{\rm c}}$, where $\mathcal H_{\text{\rm at}}$ and $\mathcal H_{\text{\rm c}}$ are 
  hyperinvariant subspaces of $U$, $\mathcal H_{\text{\rm at}}=\oplus_j \mathcal H_j$, $U|_{\mathcal H_{\text{\rm at}}}=\oplus_j \lambda_jI_{\mathcal H_j}$, 
and $U|_{\mathcal H_{\text{\rm c}}}$ has no eigenvalues. Moreover, there exist a subsequence $\{n_k\}_k$ of positive integers and a family 
 $\{\xi_j\}_j\subset \mathbb T$ such that 
\begin{equation*}
\begin{aligned}
& U^{n_k}x\to_k(\oplus_j \xi_jI_{\mathcal H_j})x\ \ \ \text{ for every }\ x \in\mathcal H_{\text{\rm at}} \\ and \ \ 
  & U^{n_k}\to_k(\oplus_j \xi_jI_{\mathcal H_j})\oplus \mathbb O_{\mathcal H_{\text{\rm c}}} \ \ \text{ weakly.}\end{aligned}\end{equation*}
\end{lemma}
\begin{proof} By Wiener theorem (see, for example, {\cite[I.7.13, Remark {\bf b}]{katz}}), there exists a  subsequence $\{n_l\}_l$ of positive integers such that  
\begin{equation}
\label{aaa1}
U^{n_l}|_{\mathcal H_{\text{\rm c}}}\to_l \mathbb O_{\mathcal H_{\text{\rm c}}}\ \text{ weakly. }
\end{equation}
There exist  a subsequence $\{l_k\}_k$ of positive integers and 
a family  $\{\xi_j\}_j\subset \mathbb T$ such that $\lambda_j^{n_{l_k}}\to_k\xi_j$ for every $j$. 
Clearly, \begin{equation}
\label{aaa2} 
U^{n_{l_k}}|_{\mathcal H_{\text{\rm at}}}x\to_k(\oplus_j \xi_jI_{\mathcal H_j})x \ \ \text{ for every } 
\ x\in\mathcal H_{\text{\rm at}}.\end{equation}
 The conclusion of the lemma follows from (\ref{aaa1}) and (\ref{aaa2}).
\end{proof}

 \begin{proposition}
\label{lll2} Let  $T\in\mathcal L(\mathcal H)$ be an invertible operator, and let $C>0$. Set 
\begin{equation*}
\begin{aligned}\mathcal X=&\{x\in\mathcal H\ :\ \text{ there exist subsequences } \ \{n_k\}_k \ 
\text{ and } \ \{l_k\}_k  \\ 
& \ \text{ of positive integers such that } \ \lim_k\|T^{-l_k}x\| \leq C\lim_k\|T^{n_k}x\|\}\end{aligned}\end{equation*}
(we suppose that the limits in the definition of $\mathcal X$ exist).
Suppose that $\operatorname{clos}\mathcal X=\mathcal H$. Then $\sup_{n\geq 0}\|T^{-n}\|\leq C (\sup_{n\geq 0}\|T^n\|)^2$.
\end{proposition}
\begin{proof} Set $M=\sup_{n\geq 0}\|T^n\|$. For $n\geq 0$, $x\in\mathcal X$, and all sufficiently large $k$ we have 
$$\|T^{-n}x\|=\|T^{l_k-n}T^{-l_k}x\|\leq M\|T^{-l_k}x\|.$$
Therefore, $$\|T^{-n}x\|\leq M \lim_k\|T^{-l_k}x\|\leq M C\lim_k\|T^{n_k}x\|\leq  CM^2\|x\|.$$  Since $\operatorname{clos}\mathcal X=\mathcal H$, 
the latest estimate takes place for every $x\in\mathcal H$.
\end{proof}

\section{ Multipliers as intertwining operators: estimates of the norms of powers of intertwined operators}

We need the following simple lemma.

\begin{lemma}
\label{lll6} Suppose $\nu$ is a positive measure, $\{\psi_n\}_n\subset L^\infty(\nu)$, $\psi\in L^\infty(\nu)$, 
$\sup_n\|\psi_n\|_\infty<\infty$, $\psi_n\to \psi$ a.e.,  $\{\eta_n\}_n\subset L^2(\nu)$, $g$, $\eta\in L^2(\nu)$, 
and $\eta_n\to \eta$ weakly in $L^2(\nu)$. Then $(g\psi_n,\eta_n)_{L^2(\nu)}\to (g\psi,\eta)_{L^2(\nu)}$.
\end{lemma}
\begin{proof} We have $$\int g\psi_n\overline\eta_n\text{\rm d}\nu -\int g\psi\overline\eta\text{\rm d}\nu  = 
\int g(\psi_n-\psi)\overline\eta_n\text{\rm d}\nu + \int g\psi(\overline\eta_n-\overline\eta)\text{\rm d}\nu.$$
Second summand tends to zero due to the weak convergence $\eta_n\to \eta$. First summand can be estimated as follows:
\begin{equation*}
\begin{aligned}\Bigl|\int g(\psi_n-\psi)\overline\eta_n\text{\rm d}\nu\Bigr| & \leq\int |g(\psi_n-\psi)\overline\eta_n|\text{\rm d}\nu \\ 
& \leq
 \Bigl(\int|g(\psi_n-\psi)|^2\text{\rm d}\nu\Bigr)^{1/2} \Bigl(\int|\eta_n|^2\text{\rm d}\nu\Bigr)^{1/2}.\end{aligned}\end{equation*}
We have $$ \sup_n\Bigl(\int|\eta_n|^2\text{\rm d}\nu\Bigr)^{1/2}<\infty,$$
because  the sequence $\{\eta_n\}_n$ converges weakly  in $L^2(\nu)$. 
First factor tends to zero by Lebesgue's dominated convergence theorem. 
\end{proof}

The following proposition is the main result of this section.

\begin{proposition}
\label{ppp1} Suppose $\mu$ is a positive Borel singular measure on $\mathbb T$, $\mu(\mathbb T)=1$, 
$\chi(\zeta)=\zeta$, $\zeta\in\overline{\mathbb D}$, $U_\mu$
is the operator of multiplication by $\chi$ in $L^2(\mu)$, $\gamma\in L^\infty(\mu)$, $W_\gamma$ is 
 the operator of multiplication by $\gamma$ in $L^2(\mu)$, and there exists a 
subsequence $\{n_k\}_k$ of positive integers such that $$U_\mu^{n_k}\to_k W_\gamma \ \ \text{ weakly.}$$
Furthermore, suppose $\theta$ is defined by $\mu$ as in (\ref{yy0}), 
$J_{\theta,1}^{-1}\in\mathcal L(L^2(\mu), \mathcal K_\theta)$ acts by the formula (\ref{yy1}),
$g\in L^2(m)$, $$X\in\mathcal L(\mathcal K_\theta, L^2(m)) \ \text{  acts by the formula } \ Xf=gf, \ \  f\in\mathcal  K_\theta.$$ 
Finally, let $p$ be an (analytic) polynomial. 
Then 
$$\lim_k \bigl( \|XJ_{\theta,1}^{-1}U_\mu^{n_k}p\|^2 -\|pXJ_{\theta,1}^{-1}\chi^{-n_k}\|^2\bigr)= \|XJ_{\theta,1}^{-1}W_\gamma p\|^2 - \|pXJ_{\theta,1}^{-1}\overline\gamma\|^2$$ 
and 
$$\lim_k \bigl(\|XJ_{\theta,1}^{-1}U_\mu^{-n_k}p\|^2-\|pXJ_{\theta,1}^{-1}\chi^{-n_k}\|^2\bigr)=\|XJ_{\theta,1}^{-1}W_\gamma^\ast p\|^2 - \|pXJ_{\theta,1}^{-1}\overline\gamma\|^2.$$
\end{proposition}
\begin{proof} Set  $\kappa_n=J_{\theta,1}^{-1}\chi^n$, $n\in\mathbb Z$. By Lemma \ref{lll7},
$$\kappa_n=\theta\overline\kappa_{-n}+ \widehat\mu(-n)(1-\theta), \ \ \ n\in\mathbb N.$$ Therefore, 
\begin{equation}
\label{yy4}\begin{aligned}
|\kappa_n|^2=&|\kappa_{-n}|^2+ 2\operatorname{Re}\bigl(\theta\overline\kappa_{-n}
\overline{\widehat\mu(-n)}(1-\overline\theta)\bigr)
 + |1-\theta|^2|\widehat\mu(-n)|^2 
\\&\ \ \text{ a.e. on } \mathbb T \ \text{(with respect to }m), \ \ \ n\in\mathbb N. \end{aligned}\end{equation}

For  $z\in\overline{\mathbb D}$ set 
$$f_z(\zeta)=\overline{\zeta\frac{p(\zeta)-p(z)}{\zeta-z}} = 
\overline{\zeta\sum_{l=1}^{\deg p} \widehat p(l)\sum_{s=0}^{l-1}\zeta^{l-1-s}z^s}, \ \ \zeta\in\mathbb T, $$ 
and for $n\in\mathbb Z$ and $z\in\mathbb D$ set 
\begin{equation}
\label{yy2}\psi_n(z)=(1-\theta(z))(\chi^n,f_z)_{L^2(\mu)}.\end{equation}
 There exists $C>0$ which depends on $p$, but does not depend on $n\in\mathbb Z$ and  $z\in\mathbb D$,
such that \begin{equation}
\label{yy11}|\psi_n(z)|\leq C \ \ \ \text{ for all } \ z\in\mathbb D, \  \ n\in\mathbb Z.\end{equation}
In particular, $\psi_n\in H^\infty$. Furthermore,  (\ref{yy2}) is fulfilled for nontangential boundary values of $\psi_n$. 
It follows  from (\ref{yy1}) that 
$$(J_{\theta,1}^{-1}U_\mu^n p)(z)=\psi_n(z)+p(z)\kappa_n(z), \ \ \ z\in\mathbb D, \ \ n\in\mathbb Z.$$ Thus, 
\begin{equation}
\label{yy5}XJ_{\theta,1}^{-1}U_\mu^n p=g\psi_n+gp\kappa_n, \ \ n\in\mathbb Z. \end{equation}
We infer from (\ref{yy4}) and (\ref{yy5}) that
\begin{equation}
\label{yy6}\begin{aligned}
\|XJ_{\theta,1}^{-1}U_\mu^n p\|^2=&\|g\psi_n\|^2+2\operatorname{Re}(g\psi_n,gp\kappa_n)+\|gp\kappa_{-n}\|^2 \\&
 + 2\operatorname{Re}\Bigl(\overline{\widehat\mu(-n)}\int_{\mathbb T}|g|^2|p|^2\theta\overline\kappa_{-n}(1-\overline\theta)\text{\rm d}m\Bigr)
 \\&+ |\widehat\mu(-n)|^2\int_{\mathbb T}|g|^2|p|^2|1-\theta|^2\text{\rm d}m,  \ \ \ n\in\mathbb N,   \end{aligned}\end{equation}
while 
\begin{equation}
\label{yy7}\|XJ_{\theta,1}^{-1}U_\mu^{-n} p\|^2=\|g\psi_{-n}\|^2+2\operatorname{Re}(g\psi_{-n},gp\kappa_{-n})+\|gp\kappa_{-n}\|^2,  \ \ \ n\in\mathbb N.  \end{equation}

Set $a=\int_{\mathbb T}\gamma\text{\rm d}\mu$. We have 
\begin{equation}
\label{yy9} \lim_k\widehat\mu(-n_k) = a. \end{equation}
Furthermore,  $\kappa_{n_k}\to_k J_{\theta,1}^{-1}\gamma$ and $\kappa_{-n_k}\to_k J_{\theta,1}^{-1}\overline\gamma$ weakly in $\mathcal K_\theta$.
 Since  $X$ is bounded,  
 \begin{equation}
\label{yy12}g\kappa_{n_k}\to_k gJ_{\theta,1}^{-1}\gamma \ \ \text{ and } \ \ g\kappa_{-n_k}\to_k gJ_{\theta,1}^{-1}\overline\gamma \ \text{  weakly in } \ L^2(m).\end{equation}
  
Since (\ref{yy2}) is fulfilled for nontangential boundary values of $\psi_n$, 
\begin{equation}
\label{yy13}\begin{aligned}\psi_{n_k}\to_k J_{\theta,1}^{-1}(p\gamma)-pJ_{\theta,1}^{-1}\gamma \ \  & \text{ and } \ \ \psi_{-n_k}\to_k J_{\theta,1}^{-1}(p\overline\gamma)-pJ_{\theta,1}^{-1}\overline\gamma \\&
 \text{ a.e. with respect to }m.\end{aligned}\end{equation}
By Lebesgue's dominated convergence theorem, \begin{equation}
\label{yy15}
\lim_k\|g\psi_{n_k}\|^2=\int_{\mathbb T}|g|^2|J_{\theta,1}^{-1}(p\gamma)-pJ_{\theta,1}^{-1}\gamma|^2\text{\rm d}m
\end{equation}
  and 
\begin{equation}
\label{yy16} \lim_k\|g\psi_{-n_k}\|^2=\int_{\mathbb T}|g|^2|J_{\theta,1}^{-1}(p\overline\gamma)-pJ_{\theta,1}^{-1}\overline\gamma|^2\text{\rm d}m.\end{equation}
By (\ref{yy11}), (\ref{yy13}) and (\ref{yy12}), Lemma \ref{lll6} can be applied to $\{\psi_{n_k}\}_k\subset L^\infty(m)$, $\{pg\kappa_{n_k}\}_k\subset L^2(m)$, 
and $g\in L^2(m)$. Thus,
\begin{equation}
\label{yy17}\begin{aligned} \lim_k(g\psi_{n_k},gp\kappa_{n_k})&=\bigl(g(J_{\theta,1}^{-1}(p\gamma)-pJ_{\theta,1}^{-1}\gamma),gpJ_{\theta,1}^{-1}\gamma\bigr) \\& =
(gJ_{\theta,1}^{-1}(p\gamma),gpJ_{\theta,1}^{-1}\gamma)-\|gpJ_{\theta,1}^{-1}\gamma\|^2. \end{aligned}\end{equation}
Similarly,  \begin{equation}
\label{yy18}\lim_k(g\psi_{-n_k},gp\kappa_{-n_k})=(gJ_{\theta,1}^{-1}(p\overline\gamma),gpJ_{\theta,1}^{-1}\overline\gamma)-\|gpJ_{\theta,1}^{-1}\overline\gamma\|^2.\end{equation}

From (\ref{yy6}), (\ref{yy9}), (\ref{yy12}),  (\ref{yy15}), and (\ref{yy17}), using elementary calculation, we conclude that 
\begin{equation}
\label{yy19}\begin{aligned}\lim_k&\Bigl(\|XJ_{\theta,1}^{-1}U_\mu^{n_k}p\|^2-\|gp\kappa_{-n_k}\|^2\Bigr) \\&
=\int_{\mathbb T}|g|^2|J_{\theta,1}^{-1}(p\gamma)|^2\text{\rm d}m-
\int_{\mathbb T}|g|^2|p|^2|J_{\theta,1}^{-1}\gamma|^2\text{\rm d}m \\&
+
2\operatorname{Re}\Bigl(\overline a\int_{\mathbb T}|g|^2|p|^2(\theta-1)
\overline{J_{\theta,1}^{-1}\overline\gamma}\text{\rm d}m\Bigr)+
|a|^2\int_{\mathbb T}|g|^2|p|^2|1-\theta|^2\text{\rm d}m,\end{aligned}\end{equation}
while from (\ref{yy7}), (\ref{yy16}), and (\ref{yy18}),  we conclude that 
\begin{equation*}
\begin{aligned}\lim_k&\Bigl(\|XJ_{\theta,1}^{-1}U_\mu^{-n_k}p\|^2-\|gp\kappa_{-n_k}\|^2\Bigr)\\&=\int_{\mathbb T}|g|^2|J_{\theta,1}^{-1}(p\overline\gamma)|^2\text{\rm d}m -
\int_{\mathbb T}|g|^2|p|^2|J_{\theta,1}^{-1}\overline\gamma|^2\text{\rm d}m.\end{aligned}\end{equation*}
The second equality from the conclusion of the proposition is proved.

Set  $u=J_{\theta,1}^{-1}(\overline\chi\overline\gamma)$. By Lemma \ref{lll7},  
$$J_{\theta,1}^{-1}\gamma = \theta \overline\chi\overline u \ \ \text{ and } \ \ J_{\theta,1}^{-1}\overline\gamma = \chi u +\overline a(1-\theta).$$
 From these equalities we conclude that 
\begin{equation}
\label{yy20}\begin{aligned}\int_{\mathbb T}&|g|^2|p|^2|J_{\theta,1}^{-1}\gamma|^2\text{\rm d}m -
2\operatorname{Re}\Bigl(\overline a\int_{\mathbb T}|g|^2|p|^2(\theta-1)
\overline{J_{\theta,1}^{-1}\overline\gamma}\text{\rm d}m\Bigr)\\& -
|a|^2\int_{\mathbb T}|g|^2|p|^2|1-\theta|^2\text{\rm d}m =\int_{\mathbb T}|g|^2|p|^2|J_{\theta,1}^{-1}\overline\gamma|^2\text{\rm d}m.\end{aligned}\end{equation}
The first equality from the conclusion of the proposition follows from (\ref{yy19}) and (\ref{yy20}).
\end{proof}

The following theorem is a corollary of Proposition \ref{ppp1}.

\begin{theorem}
\label{ttt1} Suppose $\mu$ is a positive Borel singular measure on $\mathbb T$, $\mu(\mathbb T)=1$, $U_\mu$
is the operator of multiplication by the independent variable in $L^2(\mu)$,
$\theta$ is defined by $\mu$ as in (\ref{yy0}), 
$J_{\theta,1}^{-1}\in\mathcal L(L^2(\mu), \mathcal K_\theta)$ acts by the formula (\ref{yy1}), $\omega\in H^\infty$ is an inner function, 
$g\in H^2$,
 $$X\in\mathcal L(\mathcal K_\theta, \mathcal K_\omega) \ \ \text{  acts by the formula } \ \ Xf=gf,  \ \ \ f\in \mathcal K_\theta,$$
and $\operatorname{clos}X\mathcal K_\theta = \mathcal K_\omega$. 
Let $T\in\mathcal L(\mathcal K_\omega)$ be invertible and such that $$XJ_{\theta,1}^{-1}U_\mu=TXJ_{\theta,1}^{-1}.$$
Then $$\sup_{n\geq 0}\|T^{-n}\|\leq (\sup_{n\geq 0}\|T^n\|)^5.$$
\end{theorem}

\begin{proof} We have $\mu=\mu_{\text{\rm at}}+\mu_{\text{\rm c}}$, where  $\mu_{\text{\rm at}}$ and $\mu_{\text{\rm c}}$ are the pure 
atomic and continuous parts of $\mu$, respectively. As in Proposition \ref{ppp1}, for $\gamma\in L^\infty(\mu)$,  denote by $W_\gamma$  
 the operator of multiplication by $\gamma$ acting in $L^2(\mu)$. 
 By Lemma \ref{lll1}, 
there exist a subsequence $\{n_k\}_k$ of positive integers and a function $\gamma \in L^\infty(\mu)$ such that $U_\mu^{n_k}\to_k W_\gamma$ weakly, 
 $ W_\gamma|_{L^2(\mu_{\text{\rm c}})}=\mathbb O_{L^2(\mu_{\text{\rm c}})}$,  $W_\gamma|_{L^2(\mu_{\text{\rm at}})}$ is unitary, and 
 $U_\mu^{n_k}x\to_kW_\gamma x$ for every $x\in L^2(\mu_{\text{\rm at}})$. 

Put  $$U_1=U_\mu|_{L^2(\mu_{\text{\rm at}})}, \ \ \ \ W=W_\gamma|_{L^2(\mu_{\text{\rm at}})}, $$
$$  
\mathcal M=\operatorname{clos}XJ_{\theta,1}^{-1}L^2(\mu_{\text{\rm at}}), \ \ \ \ X_1=XJ_{\theta,1}^{-1}|_{L^2(\mu_{\text{\rm at}})}.$$ 
We have $T\mathcal M\subset\mathcal M$. Put $T_1=T|_{\mathcal M}$. We have that 
$X_1\in \mathcal L(L^2(\mu_{\text{\rm at}}), \mathcal M)$ is a quasiaffinity, and $X_1U_1=T_1X_1$. 
Let $R\in\mathcal L( \mathcal M)$ be an invertible operator from  Lemma \ref{lll3} applied to  $U_1$, $W$, $T_1$ and $X_1$. We have $RX_1=X_1W$.

Suppose $p$ is a polynomial. Set $$F(p)=-\|pXJ_{\theta,1}^{-1}\overline\gamma\|^2.$$
We have $T^nXJ_{\theta,1}^{-1}p=XJ_{\theta,1}^{-1}U_\mu^np$ for all $n\in\mathbb Z$. By Proposition \ref{ppp1},  
\begin{equation*}
\begin{aligned}\lim_k & (\|T^{n_k}XJ_{\theta,1}^{-1}p\|^2-\|pXJ_{\theta,1}^{-1}\chi^{-n_k}\|^2)\\ & = \lim_k(\|XJ_{\theta,1}^{-1}U_\mu^{n_k}p\|^2-\|pXJ_{\theta,1}^{-1}\chi^{-n_k}\|^2)=\|XJ_{\theta,1}^{-1}W_\gamma p\|^2 +F(p)\\ &=\|XJ_{\theta,1}^{-1}W_\gamma P_{L^2(\mu_{\text{\rm at}})}p\|^2 +F(p)=\|XJ_{\theta,1}^{-1}WP_{L^2(\mu_{\text{\rm at}})} p\|^2 
+F(p)\\ &
=\|X_1WP_{L^2(\mu_{\text{\rm at}})} p\|^2+F(p)=\|RX_1P_{L^2(\mu_{\text{\rm at}})} p\|^2+F(p).\end{aligned}\end{equation*}
Similarly, 
$$\lim_k(\|T^{-n_k}XJ_{\theta,1}^{-1}p\|^2-\|pXJ_{\theta,1}^{-1}\chi^{-n_k}\|^2) =  \|R^{-1}X_1 P_{L^2(\mu_{\text{\rm at}})} p\|^2+F(p).$$

Set $M=\sup_{n\geq 0}\|T^n\|$. Since  $\sup_k\|pXJ_{\theta,1}^{-1}\chi^{-n_k}\|<\infty$, there exists a subsequence 
$\{k_j\}_j$ of positive integers (which can depend on $p$) such that $\lim_j\|pXJ_{\theta,1}^{-1}\chi^{-n_{k_j}}\|$
exists. Set $$L(p)=\lim_j\|pXJ_{\theta,1}^{-1}\chi^{-n_{k_j}}\|^2.$$ Since $\chi^{-n_{k_j}}\to\overline\gamma$ weakly, $L(p)+F(p)\geq 0$. We have
\begin{equation}\label{qqq1}\begin{aligned}\lim_j\|T^{n_{k_j}}XJ_{\theta,1}^{-1}p\|^2&=\lim_j(\|T^{n_{k_j}}XJ_{\theta,1}^{-1}p\|^2-\|pXJ_{\theta,1}^{-1}\chi^{-n_{k_j}}\|^2) + L(p) \\&=\|RX_1P_{L^2(\mu_{\text{\rm at}})} p\|^2+F(p)+ L(p)\end{aligned}\end{equation} 
and
\begin{equation}\label{qqq2}\begin{aligned}
\lim_j\|T^{-n_{k_j}}XJ_{\theta,1}^{-1}p\|^2&=\lim_j(\|T^{-n_{k_j}}XJ_{\theta,1}^{-1}p\|^2-\|pXJ_{\theta,1}^{-1}\chi^{-n_{k_j}}\|^2) + L(p)\\&
=
\|R^{-1}X_1 P_{L^2(\mu_{\text{\rm at}})} p\|^2+F(p) + L(p).\end{aligned}\end{equation} 
By Lemma \ref{lll3}, $\|R^{-1}X_1 P_{L^2(\mu_{\text{\rm at}})} p\|\leq M^3\|RX_1P_{L^2(\mu_{\text{\rm at}})} p\|$.
From this estimate, (\ref{qqq1}) and (\ref{qqq2}) we conclude that $$\lim_j\|T^{-n_{k_j}}XJ_{\theta,1}^{-1}p\|^2\leq M^6\lim_j\|T^{n_{k_j}}XJ_{\theta,1}^{-1}p\|^2.$$

By Proposition \ref{lll2} applied with $\mathcal X=\{XJ_{\theta,1}^{-1}p\ :\ p \text{ is a polynomial}\}$ and $C=M^3$, 
$$\sup_{n\geq 0}\|T^{-n}\|\leq M^3(\sup_{n\geq 0}\|T^n\|)^2=(\sup_{n\geq 0}\|T^n\|)^5.$$
\end{proof}

\section{Asymmetric truncated Toeplitz operators as intertwining operators}

In this section, for inner functions $\theta$ and $\omega$, and $g\in L^2(m)$,  operators acting from $\mathcal K_\theta$ to $\mathcal K_\omega$ 
 by the formula $f\mapsto P_{\mathcal K_\omega}gf$,  $f\in \mathcal K_\theta$, are considered. These operators are called {\it 
asymmetric truncated Toeplitz operators}, see \cite{cp}. If $\theta=\omega$, then these operators are called {\it 
 truncated Toeplitz operators}, see, for example,  \cite{cft}, \cite{gr},  {\cite[Sec. 12]{grm}} and references therein. Recall that  $U_{(\theta)1}$ is  defined in ({\ref{y3}). It will be shown that if an asymmetric truncated Toeplitz operator  intertwines $U_{(\theta)1}$ with a power bounded operator $T$ and  some additional conditions are fulfilled, then the estimate of the norms of negative 
powers of $T$ depends on the estimate of the norms of positive powers of $T$ only. Author do not know 
if the additional conditions can be dropped. 

\begin{theorem}\label{ttt2}  Suppose  $\theta$, $\omega\in H^\infty$ are inner function, $\theta(0)=0$, and $U_{(\theta)1}$ has no eigenvalues.  Furthermore, suppose 
$g\in L^2(m)$, $X\in\mathcal L(\mathcal K_\theta, \mathcal K_\omega)$ acts by the formula 
$$Xf=P_{\mathcal K_\omega}gf, \ \ \ \ f\in \mathcal K_\theta,$$
$\operatorname{clos}X\mathcal K_\theta = \mathcal K_\omega$, $K\in\mathcal L(\mathcal K_\theta, L^2(m)\ominus\mathcal K_\omega)$ acts by the formula 
$$Kf=P_{L^2(m)\ominus\mathcal K_\omega}gf, \ \ \ f\in \mathcal K_\theta,$$ and $K$ is compact.
Let $T\in\mathcal L(\mathcal K_\omega)$ be invertible and such that $$XU_{(\theta)1}=TX.$$
Then $$\sup_{n\geq 0}\|T^{-n}\|\leq (\sup_{n\geq 0}\|T^n\|)^2.$$
\end{theorem}

\begin{proof}  Proposition \ref{ppp1} will be applied. For convenience,  denote by $\mu$ a positive Borel singular measure on $\mathbb T$ which is defined by $\theta$ in (\ref{yy0}). 
By (\ref{y5}),  $U_\mu\cong U_{(\theta)1}$. Therefore, $U_\mu$ has no eigenvalues. By Wiener theorem (see, for example, {\cite[I.7.13, Remark {\bf b}]{katz}}), there exists a subsequence $\{n_k\}_k$ of positive integers such that 
 $U_\mu^{n_k}\to_k \mathbb O$ weakly. 
Define $X_1\in\mathcal L(\mathcal K_\theta, L^2(m))$ by the formula 
$$X_1 f= Xf + Kf=gf, \ \ \ f\in \mathcal K_\theta.$$ Recall that $\chi(\zeta)=\zeta$, $\zeta\in\overline{\mathbb D}$. By Proposition \ref{ppp1}, 
\begin{equation}\label{qqq4}\begin{aligned}&\lim_k \bigl( \|X_1J_{\theta,1}^{-1}U_\mu^{n_k}p\|^2 -\|pX_1J_{\theta,1}^{-1}\chi^{-n_k}\|^2\bigr)= 0 \\
\text{and }\ \ &\lim_k \bigl(\|X_1J_{\theta,1}^{-1}U_\mu^{-n_k}p\|^2-\|pX_1J_{\theta,1}^{-1}\chi^{-n_k}\|^2\bigr)=0\end{aligned}\end{equation} 
for every (analytic) polynomial $p$.

Since  $K$ is compact, $\lim_k\|KJ_{\theta,1}^{-1}U_\mu^{n_k}x\|=0$ and $\lim_k\|KJ_{\theta,1}^{-1}U_\mu^{-n_k}x\|=0$  for every $x\in L^2(\mu)$.  
Therefore, (\ref{qqq4}) can be rewritten as
\begin{equation}\label{qqq5}\begin{aligned}&\lim_k \bigl( \|XJ_{\theta,1}^{-1}U_\mu^{n_k}p\|^2 -\|pXJ_{\theta,1}^{-1}\chi^{-n_k}\|^2\bigr)= 0 \\ \text{and }\ \ &\lim_k \bigl(\|XJ_{\theta,1}^{-1}U_\mu^{-n_k}p\|^2-\|pXJ_{\theta,1}^{-1}\chi^{-n_k}\|^2\bigr)=0 \end{aligned}\end{equation}
for every (analytic) polynomial $p$.

 Suppose $p$ is a polynomial. Since  $\sup_k\|pXJ_{\theta,1}^{-1}\chi^{-n_k}\|<\infty$, there exists a subsequence 
$\{k_j\}_j$ of positive integers (which can depend on $p$) such that  $\lim_j\|pXJ_{\theta,1}^{-1}\chi^{-n_{k_j}}\|$
exists. By (\ref{y5}) and (\ref{qqq5}), 
$$\lim_j\|T^{n_{k_j}}XJ_{\theta,1}^{-1}p\|=\lim_j\|XJ_{\theta,1}^{-1}U_\mu^{n_{k_j}}p\|=\lim_j\|pXJ_{\theta,1}^{-1}\chi^{-n_{k_j}}\|.$$ Similarly,  
$$\lim_j\|T^{-n_{k_j}}XJ_{\theta,1}^{-1}p\|=\lim_j\|pXJ_{\theta,1}^{-1}\chi^{-n_{k_j}}\|.$$

By Proposition \ref{lll2} applied with $\mathcal X=\{XJ_{\theta,1}^{-1}p\ :\ p \text{ is a polynomial}\}$ and $C=1$, 
$$\sup_{n\geq 0}\|T^{-n}\|\leq (\sup_{n\geq 0}\|T^n\|)^2.$$
\end{proof}

The following lemma gives  sufficient conditions for the operator 
\begin{equation*}\begin{aligned}&K\in\mathcal L(\mathcal K_\theta, L^2(m)\ominus\mathcal K_\omega) \ \text{ acting by the formula } \\&
Kf=P_{L^2(m)\ominus\mathcal K_\omega}gf, \ \  \  f\in \mathcal K_\theta, \ \  \ g\in L^2(m),\end{aligned}\end{equation*} be  compact. Recall that 
$$QC=(H^\infty+C(\mathbb T))\cap(\overline{H^\infty+C(\mathbb T)})$$ is the algebra of {\it quasicontinuous} functions, see, for example, {\cite[Appendix 4, \S 53]{n}}, {\cite[\S 3.3]{p}}.

\begin{lemma}
\label{lll9} Suppose $\theta$, $\omega\in H^\infty$ are inner function, $\omega\overline\theta\in H^\infty+C(\mathbb T)$, and $g\in QC$.
Then the operator  $K\in\mathcal L(\mathcal K_\theta, L^2(m)\ominus\mathcal K_\omega)$ acting by the formula 
$Kf=P_{L^2(m)\ominus\mathcal K_\omega}gf$,  $f\in \mathcal K_\theta$, is compact.
\end{lemma}
\begin{proof} We have $ L^2(m)\ominus\mathcal K_\omega=\omega H^2\oplus H^2_-$. Set 
$$K_+=P_{\omega H^2}K\in\mathcal L(\mathcal K_\theta,\omega H^2) \ \ \text{ and } \ \ \ K_-= P_{H^2_-}K\in\mathcal L(\mathcal K_\theta, H^2_-).$$ 
Then $K$ is compact if and only if $K_+$ and $K_-$ are compact.

Define $K_{+\ast}\in\mathcal L(\omega H^2, \theta H^2_-)$ by the formula
$$ K_{+\ast}(\omega h) = P_{\theta H^2_-}\overline g \omega h, \ \ \ h\in H^2.$$
Since $P_{\theta H^2_-}f=\theta P_{H^2_-}\overline\theta f$ for every $f\in L^2(m)$, we have 
$$ K_{+\ast}(\omega h) = \theta P_{H^2_-}\overline\theta\overline g \omega h, \ \ \ h\in H^2.$$
By Hartman's theorem (see, for example, {\cite[Appendix 4, \S 22]{n}}, {\cite[Theorem 1.5.5]{p}}), $K_{+\ast}$ is compact if and only if 
\begin{equation}
\label{iii1}\overline\theta\overline g \omega\in H^\infty+C(\mathbb T).
\end{equation} Since $H^\infty+C(\mathbb T)$ is an algebra, (\ref{iii1}) follows from the conditions of the lemma. Clearly,  $K_+= (K_{+\ast})^\ast|_{\mathcal K_\theta}$, therefore, $K_+$ is compact.

Define $K_1\in\mathcal L(H^2, H^2_-)$ by the formula 
$$ K_1 h = P_{H^2_-} g h, \ \ \ h\in H^2.$$
By Hartman's theorem (see, for example, {\cite[Appendix 4, \S 22]{n}}, {\cite[Theorem 1.5.5]{p}}), $K_1$ is compact if and only if $ g \in H^\infty+C(\mathbb T)$. 
Since $K_-=K_1|_{\mathcal K_\theta}$, $K_-$ is compact.
\end{proof}

In {\cite[Sec. 8]{cft}}, some results on  invertibility of truncated Toeplitz operators are given. Recall the simplest result. 
Let $\theta$ be an inner function, and let $g\in H^\infty$. Let  $X$ be  truncated Toeplitz operator with the symbol $g$. That is, 
 $X\in\mathcal L(\mathcal K_\theta)$, and $X$  acts  by the formula
$Xf=P_{\mathcal K_\theta}gf$, $f\in\mathcal K_\theta$. Denote by $\alpha$ the greatest common divisor of the inner factor of $g$ and $\theta$. 
We have  $\ker X=\overline\alpha\theta\mathcal K_\alpha$ and 
$\operatorname{clos}X \mathcal K_\theta = \alpha\mathcal K_{\overline\alpha\theta}$. 
Furthermore,  if $1/g\in H^\infty$, then $X$ is invertible and $X^{-1}f=P_{\mathcal K_\theta}((1/g)f)$, $f\in\mathcal K_\theta$. 
The following lemma is a simple generalization of these facts to the case of asymmetric truncated Toeplitz operators.

\begin{lemma}
\label{lll31} Let $\theta$ and $\omega$ be inner functions, and let $g\in H^\infty$. Let $X$ be asymmetric truncated Toeplitz operator with the symbol $g$. That is,   $X\in\mathcal L(\mathcal K_\theta,\mathcal K_\omega)$, and $X$ acts 
by the formula $Xf=P_{\mathcal K_\omega}gf$, $f\in\mathcal K_\theta$. Denote by $\alpha$ the greatest common divisor of the inner factor of $g$ 
and $\omega$. Then $\ker X=\mathcal K_\theta\cap\overline\alpha\omega H^2$, and if $\mathcal K_\theta\vee\overline\alpha\omega H^2=H^2$, then 
 $\operatorname{clos}X \mathcal K_\theta = \alpha\mathcal K_{\overline\alpha\omega}$. 

Furthermore, suppose  $1/g\in H^\infty$. Then $X \mathcal K_\theta = \mathcal K_\omega$ if and only if $H^2=\mathcal K_\theta + \omega H^2$, and  $X$ is invertible if and only if \begin{equation} \label{rrrtoep} H^2=\mathcal K_\theta\dotplus \omega H^2. \end{equation} Moreover,  
$X^{-1}f=\mathcal P_{\mathcal K_\theta|| \omega H^2}((1/g)f)$, $f\in\mathcal K_\omega$, where $\mathcal P_{\mathcal K_\theta|| \omega H^2}$ is the (non-orthogonal) 
projection on $\mathcal K_\theta$ parallel to $\omega H^2$. 
\end{lemma}
\begin{proof} Suppose $f\in\ker X$. It means that $gf\in\omega H^2$. Since the inner factor of $\overline\alpha g$ and $\overline\alpha\omega$
 are relatively prime, we conclude that  $f\in\overline\alpha\omega H^2$. Conversely, if $f\in\mathcal K_\theta\cap\overline\alpha\omega H^2$,
then $gf\in\omega H^2$, therefore, $Xf=0$. 

Let $g_e$ be the outer factor of $g$. Then $g=\alpha\beta g_e$, where $\beta$ is an inner function relatively prime with 
$\overline\alpha\omega$. Clearly, $P_{\mathcal K_\omega}gf=\alpha P_{\mathcal K_{\overline\alpha\omega}}\beta g_e f$ for every $f\in H^2$, 
and $P_{\mathcal K_{\overline\alpha\omega}}\beta g_e \overline\alpha\omega h=0$ for every $h\in H^2$. Define 
$X_1\in\mathcal L(H^2,\mathcal K_{\overline\alpha\omega})$ by the formula
$X_1h=P_{\mathcal K_{\overline\alpha\omega}}\beta g_e h$, $h\in H^2$. We have $X_1\overline\alpha\omega H^2=\{0\}$. 
Therefore, $$\operatorname{clos}X \mathcal K_\theta = \alpha\operatorname{clos}X_1(\mathcal K_\theta\vee\overline\alpha\omega H^2).$$ 
If $\mathcal K_\theta\vee\overline\alpha\omega H^2=H^2$, then 
$$\operatorname{clos}X \mathcal K_\theta = \alpha\operatorname{clos}X_1 H^2=
\alpha\operatorname{clos}P_{\mathcal K_{\overline\alpha\omega}}\beta H^2=\alpha\mathcal K_{\overline\alpha\omega}.$$

Suppose  that $1/g\in H^\infty$. We have 
$$H^2 =\frac{1}{g}H^2=\frac{1}{g}(\mathcal K_\omega+\omega H^2)=\frac{1}{g}\mathcal K_\omega+\omega H^2.$$ 
That is, an arbitrary function $h\in H^2$ has the form $h=(1/g)\psi + \omega h_1$, where $\psi\in\mathcal K_\omega$ and $h_1\in H^2$. 
Suppose  that  $X\mathcal K_\theta = \mathcal K_\omega$. Then there exists $f\in\mathcal K_\theta$ such that $\psi=Xf = gf+\omega h_2$ for some $h_2\in H^2$. We obtain that
$$h=\frac{1}{g}(gf+\omega h_2)+ \omega h_1=f+\omega(\frac{1}{g}h_2 + h_1).$$ It means that $H^2=\mathcal K_\theta + \omega H^2$.

Conversely, suppose that $1/g\in H^\infty$ and $H^2=\mathcal K_\theta + \omega H^2$. Let $\psi\in\mathcal K_\omega$. There exist $f\in\mathcal K_\theta$ 
and $h\in H^2$ such that $(1/g)\psi=f+\omega h$. Clearly, $\psi=Xf$. We obtain that $X \mathcal K_\theta = \mathcal K_\omega$.

If $1/g\in H^\infty$, then $g$ is outer, therefore, $\alpha$ is a constant function. Thus, $\ker X= \mathcal K_\theta\cap\omega H^2$.
Since $X$ is invertible if and only if $\ker X=\{0\}$ and $X \mathcal K_\theta = \mathcal K_\omega$, we obtain that $X$ is invertible if and only if  (\ref{rrrtoep}) is fulfilled. The formula for $X^{-1}$ can be checked straightforward.
\end{proof}

\begin{remark} Let $\theta$ and $\omega$ be inner functions, and let $T_{\overline\theta\omega}\in\mathcal L( H^2)$ be the 
 Toeplitz operator with the symbol $\overline\theta\omega$, i.e., $T_{\overline\theta\omega}h=P_{H^2}\overline\theta\omega h$, $h\in H^2$.
(On Toeplitz operators, see, for example, \cite{n}, \cite{p}.)  It is easy to see that the conditions 
 $\mathcal K_\theta\cap\omega H^2=\{0\}$ and  $H^2 = \mathcal K_\theta+\omega H^2$ are equivalent to the conditions 
$\ker T_{\overline\theta\omega}=\{0\}$ and  $H^2=T_{\overline\theta\omega} H^2$, respectively. Thus,   (\ref{rrrtoep}) is 
equvalent to invertibility of $T_{\overline\theta\omega}$. 
\end{remark}

 \begin{remark}The necessary and sufficient conditions on  inner functions  
$\theta$ and $\omega$ under which (\ref{rrrtoep}) is fulfilled  can be found in  {\cite[\S VIII.6]{n}}.
In particular, the condition $ \|\theta-\omega\|_\infty<1 $ is sufficient for  (\ref{rrrtoep}). 
If $\theta$ is a Blaschke product, that it is possible to find  a Blaschke product  $\omega$ with zeros so closed to zeros of $\theta$ 
that  $\|\theta-\omega\|_\infty<1$ and $\omega\overline\theta\in C(\mathbb T)$. Let $g\in H^\infty\cap QC$ be such that $1/g\in H^\infty$. 
Then  $\theta$, $\omega$ and $g$ satisfy the conditions of Lemma \ref{lll9}, and $X$ defined as in Lemma \ref{lll31} is invertible.
\end{remark}

\section{ Reducing one dimensional perturbations of unitaries to a special case}

\begin{lemma}
\label{lllsp} Suppose  $U\in\mathcal L(\mathcal K)$, $T$, $V$, $K\in\mathcal L(\mathcal H)$,  $X\in\mathcal L(\mathcal K,\mathcal H)$,
 $U$, $V$ are unitary, $K$ is compact, $T=V+K$, $XU=TX$, 
and $\operatorname{clos}X\mathcal K =\mathcal H$. Then $\sigma(T)\subset\mathbb T$.
\end{lemma}
\begin{proof} Since $T=V+K$ and $K$ is compact, $\sigma_e(T)=\sigma_e(V)\subset\mathbb T$. If 
$\lambda\in\sigma(T)\setminus\sigma_e(T)$, then $\operatorname{ind}(T-\lambda)=\operatorname{ind}(V-\lambda)=0$,
where $\operatorname{ind}$ is the Fredholm index. In particular, $\overline\lambda$ is an eigenvalue of $T^\ast$. Since $\ker X^\ast=\{0\}$, 
$\overline\lambda$ is an eigenvalue of $U^\ast$. Thus, $\overline\lambda\in\mathbb T$.
\end{proof}

We need some notions and facts from \cite{mv} and \cite{s}.
An operator $T$ is called {\it almost unitary}, if there exists a unitary operator $V$ such that $T-V$ is of trace class and 
$\mathbb D\not\subset \sigma(T)$. Clearly, $T$ is almost unitary if and only if $T^\ast$ is almost unitary. 
For an almost unitary operator $T\in\mathcal L(\mathcal H)$ denote by $\mathcal H_s(T)$ the set of  $x\in\mathcal H$ such that for every $y\in\mathcal H$ and for a.e. $\zeta\in\mathbb T$ with respect to $m$,
nontangential boundary values of $((T-\lambda)^{-1}x,y) $ taking inside of $\mathbb D$ and outside of $\mathbb D$ when $\lambda\to\zeta$ coincide. The space $\mathcal H_s(T)$ is called the {\it singular} subspaces of $T$.
It is easy to see that $T\mathcal H_s(T)\subset\mathcal H_s(T)$, and  it is proved in \cite{mv} that $\mathcal H_s(T)$ is closed. The space  
$\mathcal H_a(T)=\mathcal H\ominus\mathcal H_s(T^\ast)$ is called the {\it absolutely continuous} subspaces of $T$. 
If $T$ is unitary, then $T$ is singular if and only if $ \mathcal H_s(T)=\mathcal H$. It  follows exactly from the definition that 
if $ \mathcal H_s(T)=\mathcal H$ for an almost unitary operator $T$, then  $ \mathcal H_s(T^\ast)=\mathcal H$.

\begin{lemma}
\label{lllsn}  Suppose  $U\in\mathcal L(\mathcal K)$, $T\in\mathcal L(\mathcal H)$,  $X\in\mathcal L(\mathcal K,\mathcal H)$,
 $U$ is singular  unitary, $T$ is almost unitary, $XU=TX$, 
and $\operatorname{clos}X\mathcal K =\mathcal H$. Then $\mathcal H_s(T)=\mathcal H$.
\end{lemma}
\begin{proof} We have $$X(U-\lambda)^{-1}=(T-\lambda)^{-1}X \ \ \text{ for every } \lambda\notin\mathbb T\cup\sigma(T).$$
Therefore, 
\begin{equation*}\begin{aligned}((T-\lambda)^{-1}Xx,y) & =(X(U-\lambda)^{-1}x,y)= ((U-\lambda)^{-1}x,X^\ast y) \\ & \text{ for every } x\in\mathcal K \ \text{ and } 
y\in\mathcal H.\end{aligned}\end{equation*} Since $U$ is singular, we conclude that $Xx\in \mathcal H_s(T)$ for every $x\in\mathcal K$. Since $ \mathcal H_s(T)$ is closed, 
we obtain that $ \mathcal H_s(T)\supset\operatorname{clos}X\mathcal K =\mathcal H$.
\end{proof}

 The following theorem is a part of {\cite[Theorem 1.3]{s}}.

\begin{theorem}
\label{ttts}\cite{s}  Suppose $T$, $V\in\mathcal L(\mathcal H)$, $T-V$ is of trace class, $V$ is unitary, and $\mathbb D\not\subset \sigma(T)$.
Put $$\mathcal D=\{x\in\mathcal H_a(T) \ : \text{ there exists } \lim_{n\to\infty} V^{-n}T^nx \}.$$
For $x\in\mathcal D$ set $\Omega x=\lim_{n\to\infty} V^{-n}T^nx$. Then $$\operatorname{clos}\mathcal D=\mathcal H_a(T) \ \ \text{ and } \ \ 
\operatorname{clos}\Omega\mathcal D=\mathcal H_a(V).$$ \end{theorem}

\begin{corollary} 
\label{cor7.4}
Suppose $T$, $V\in\mathcal L(\mathcal H)$, $T-V$ is of trace class,  $V$ is unitary, $\mathbb D\not\subset \sigma(T)$, and 
$ \mathcal H_s(T)=\mathcal H$. Then  $V$ is a singular unitary operator. 
\end{corollary}
\begin{proof} By Theorem \ref{ttts}, $\mathcal H_a(V)=\operatorname{clos}\Omega\mathcal D$, where $\Omega$ is some linear mapping 
and $\mathcal D\subset\mathcal H_a(T)$. Since $\mathcal H_a(T)=\{0\}$, we conclude that $\mathcal H_a(V)=\{0\}$.  Consequently, $ \mathcal H_s(V^\ast)=\mathcal H$.
Therefore, $V^\ast$ is a singular  unitary operator. Indeed, $V$ is a singular  unitary operator, too.
\end{proof}

Recall that, for a positive measure $\mu$,   by $U_\mu$ the  operator of multiplication by the independent variable in $L^2(\mu)$ is denoted.
The following elementary lemma is given for a convenience of references. 

\begin{lemma}
\label{lll11} Suppose $\chi(z)=z$, $z\in\mathbb T$, $\mu$ is a positive Borel measure on $\mathbb T$, $\varphi\in L^2(\mu)$, 
and $\varphi \neq 0$ a.e. with respect to $\mu$. 
Set $a=\int|\varphi|^2\text{\rm d}\mu$ and $\text{\rm d}\mu_1=\frac{1}{a}|\varphi|^2\text{\rm d}\mu$. Define $Z\in\mathcal L(L^2(\mu),L^2(\mu_1))$ by 
the formula $Zf=\sqrt a\frac{\overline\chi}{\varphi}f$, $f\in L^2(\mu)$. Then $\mu_1(\mathbb T)=1$, $\mu$ and $\mu_1$ are mutually 
absolutely continuous, $Z$ is unitary, $ZU_\mu=U_{\mu_1}Z$ and $Z\varphi=\sqrt a\overline\chi$.
\end{lemma}

\begin{proposition}
\label{ppp01} Suppose $\chi(z)=z$, $z\in\mathbb T$, $\mu$, $\nu$ are positive  Borel measures on $\mathbb T$, 
$\varphi$, $\psi\in L^2(\nu)$,  $$T\in\mathcal L(L^2(\nu)), \ \  T=U_\nu+(\cdot, \psi)\varphi, \ \ Y\in\mathcal L(L^2(\mu),L^2(\nu)), \ \ YU_\mu=TY,$$ 
 $\psi\neq 0$ a.e. with respect to $\nu$, and $Y^\ast\psi\neq 0$ a.e. with respect to $\mu$.
Then there exist positive  Borel measures $\mu_1$, $\nu_1$ on $\mathbb T$, $\varphi_1\in L^2(\nu)$, and 
$X\in\mathcal L(L^2(\mu_1),L^2(\nu_1))$ such that $\mu_1(\mathbb T)=1$, $\nu_1(\mathbb T)=1$, $\mu$ and $\mu_1$ are mutually absolutely continuous, 
$\nu$ and $\nu_1$ are mutually absolutely continuous, $$X^\ast\overline\chi= \overline\chi, \ \  
XU_{\mu_1}=T_1X, \ \text{and } \  T\cong T_1, \ \ \text{ where } T_1=U_{\nu_1}+(\cdot, \overline\chi)\varphi_1.$$ Moreover, if 
$\varphi\neq 0$ a.e. with respect to $\nu$, then $\varphi_1\neq 0$ a.e. with respect to $\nu_1$, and if $Y$ is a quasiaffinity, then $X$ 
is a quasiaffinity.
\end{proposition}
\begin{proof} Define $\nu_1$, $a_1$ and $Z_1\in\mathcal L(L^2(\nu),L^2(\nu_1))$ as in Lemma \ref{lll11} applied to $\nu$ and $\psi$. Set 
$\varphi_1=a_1\frac{\overline\chi}{\psi}\varphi$. Then $Z_1 T=T_1Z_1$. We have $Y^\ast Z_1^{-1}\overline\chi = \frac{1}{\sqrt{a_1}}Y^\ast\psi$. 
 Define $\mu_1$, $a_2$
 and $Z_2$ as in Lemma \ref{lll11} applied to $\mu$ and $\frac{1}{\sqrt{a_1}}Y^\ast\psi$. Set $X_1=Z_1YZ_2^{-1}$. Then $X_1U_{\mu_1}=T_1X_1$ and 
$X_1^\ast\overline\chi=Z_2Y^\ast Z_1^{-1}\overline\chi =Z_2\frac{1}{\sqrt{a_1}}Y^\ast\psi=\sqrt{a_2}\overline\chi$. Set $X=\frac{1}{\sqrt{a_2}}X_1$.
Then $X$ satisfies the conclusion of the proposition.
\end{proof}

The following lemma is very simple, therefore, its proof is omitted.

\begin{lemma}
\label{lll01}   Suppose $R\in\mathcal L(\mathcal H)$, $u$, $v\in \mathcal H$, and  $T=R+(\cdot,v)u$. Set $\mathcal M=\bigvee_{n\geq 0}R^nu$. Then  
$$\mathcal M=\bigvee_{n\geq 0}T^nu, \ \ \ T|_{\mathcal M}=R|_{\mathcal M}+(\cdot,P_{\mathcal M}v)u  \ \  \text{ and } 
 \ \ P_{\mathcal M^\perp}T|_{\mathcal M^\perp}=P_{\mathcal M^\perp}R|_{\mathcal M^\perp}.$$
\end{lemma}

\begin{corollary} 
\label{ccc1} Suppose $R\in\mathcal L(\mathcal H)$, $u$, $v\in \mathcal H$, $T=R+(\cdot,v)u$, $V\in\mathcal L(\mathcal K)$, $Y\in\mathcal L(\mathcal K,\mathcal H)$,
$\ker Y=\{0\}$, and $TY=YV$. If $v$ is a cyclic vector for $R^\ast$, then $Y^\ast v$ a cyclic vector for $V^\ast$.
\end{corollary}
\begin{proof} Clearly, $T^\ast=R^\ast+(\cdot,u)v$. By Lemma \ref{lll01} applied to $R^\ast$ and $T^\ast$,   $v$ is a cyclic vector for $T^\ast$. 
Since $Y^\ast T^\ast=V^\ast Y^\ast$ and $\operatorname{clos}Y^\ast\mathcal H=\mathcal K$, we conclude that $Y^\ast v$ a cyclic vector for $V^\ast$. 
\end{proof}

\begin{lemma}
\label{lll02} Suppose $T_1\in\mathcal L(\mathcal H_1)$ and $T_2\in\mathcal L(\mathcal H_2)$ are invertible, and $T\in\mathcal L(\mathcal H_1\oplus\mathcal H_2)$ has the 
form $$T=\begin{pmatrix} T_1 & \ast\\ \mathbb O & T_2\end{pmatrix}.$$ Then $$\|T^{-n}\|\leq 
\max(1,\|T_1^{-n}\|)\cdot\max(1,\|T_2^{-n}\|)\cdot\bigl(\max(2,2\|T^n\|^2+1)\bigr)^{1/2}\ \ \text{ for } n\geq 0.$$
\end{lemma}
\begin{proof} Set $A_n=P_{\mathcal H_1}T^n|_{\mathcal H_2}$, $n\geq 1$. Then $\|A_n\|\leq\|T^n\|$ and 
$$T^{-n}=\begin{pmatrix} T_1^{-n} & -T_1^{-n}A_nT_2^{-n}\\ \mathbb O & T_2^{-n}\end{pmatrix}=
\begin{pmatrix} T_1^{-n} & \mathbb O\\ \mathbb O & I\end{pmatrix}\cdot\begin{pmatrix} I & -A_n\\ \mathbb O & I\end{pmatrix}\cdot
\begin{pmatrix} I & \mathbb O\\ \mathbb O & T_2^{-n}\end{pmatrix}.$$
The conclusion of the lemma follows from the estimate
$$\Biggl\|\begin{pmatrix} I & -A_n\\ \mathbb O & I\end{pmatrix}\Biggr\|^2 \leq \max(2,2\|A_n\|^2+1)\leq \max(2,2\|T^n\|^2+1)$$
and the equalities $$\| T_1^{-n}\oplus I\| = \max(1,\|T_1^{-n}\|) \ \ \text{  and } \ \ \|I\oplus T_2^{-n}\| = \max(1,\|T_2^{-n}\|).$$
\end{proof}

\begin{theorem}
\label{ttt5} Suppose $V\in\mathcal L(\mathcal H)$ is a reductive unitary operator, $u$, $v\in \mathcal H$, and  $$T=V+(\cdot,v)u.$$ Then there exist a 
positive Borel measure $\nu$ on $\mathbb T$, absolutely continuous with respect to the spectral measure of $V$, $\varphi$, $\psi\in L^2(\nu)$,
such that $\varphi\neq 0$ and $\psi\neq 0$ a.e. with respect to $\nu$, and  unitary operators $V_1$ and $V_2$ such that  
\begin{equation}
\label{qqq21}T\cong\begin{pmatrix}  V_1 & \ast & \ast \\ \mathbb O & T_1 & \ast\\ \mathbb O &\mathbb O & V_2\end{pmatrix}, \ \ \text{ where }
\ T_1=U_\nu+(\cdot,\psi)\varphi.
\end{equation}
Consequently, if $T$ is invertible, then 
\begin{equation}
\label{qqq22}\sup_{n\geq 0}\|T^{-n}\|\leq(\sup_{n\geq 0}\|T_1^{-n}\|)\cdot(2\sup_{n\geq 0}\|T^n\|^2+1).
\end{equation}
\end{theorem}
\begin{proof} Set $\mathcal M=\bigvee_{n\geq 0}V^nu$. Since $V$ is reductive and $\mathcal M$ is an invariant subspace of $V$, $\mathcal M^\perp$ is an 
invariant subspace of $V$, too. Therefore, $P_{\mathcal M^\perp}V|_{\mathcal M^\perp}=V|_{\mathcal M^\perp}$ is unitary. Set $V_2=V|_{\mathcal M^\perp}$.
Since $V|_{\mathcal M}$ is cyclic, 
there exist a positive Borel measure $\nu_1$ on $\mathbb T$, absolutely continuous 
with respect to the spectral measure of $V$, such that $V|_{\mathcal M}\cong U_{\nu_1}$. Let $\varphi_1\in L^2(\nu_1)$ be the function 
which correspondents to $u$ under this unitarily equivalence. Since $u$ is a cyclic vector for $V|_{\mathcal M}$,  $\varphi_1$ is a cyclic 
vector for $U_{\nu_1}$. Therefore, $\varphi_1\neq 0$ a.e. with respect to $\nu_1$. By Lemma \ref{lll01}, there exists $\psi_1\in L^2(\nu_1)$ 
such that $$T\cong\begin{pmatrix} U_{\nu_1}+ (\cdot,\psi_1)\varphi_1 & \ast\\ \mathbb O & V_2\end{pmatrix}.$$

There exists a Borel set $\tau\subset\mathbb T$ such that $\psi_1\neq 0$ a.e. on $\tau$ and $\psi_1=0$ a.e. on $\mathbb T\setminus\tau$. Set 
$$\nu=\nu_1|_\tau, \ \ \psi=\psi_1|_\tau, \ \ \varphi=\varphi_1|_\tau, \ \ V_1=U_{\nu_1}|_{L^2(\mathbb T\setminus\tau,\nu_1)}.$$
It is easy to see that  $$ U_{\nu_1}+ (\cdot,\psi_1)\varphi_1=\begin{pmatrix} V_1 & \ast\\ \mathbb O & T_1\end{pmatrix}$$
with respest to the decomposition $L^2(\nu_1)=L^2(\mathbb T\setminus\tau,\nu_1)\oplus L^2(\tau,\nu_1)$.
Thus, (\ref{qqq21}) is proved.

Suppose that $T$ is invertible. Since $V_1$ and $V_2$ are invertible,  it follows from the matrix form of $T$ that $T_1$ is invertible.
  Set $$T_2=\begin{pmatrix} T_1 & \ast\\ \mathbb O & V_2\end{pmatrix}.$$ Clearly, $\|T_2^n\|\leq\|T^n\|$ for all $n\geq 0$.
By Lemma \ref{lll02}, $$\sup_{n\geq 0}\|T_2^{-n}\|\leq (\sup_{n\geq 0}\|T_1^{-n}\|)\cdot(2\sup_{n\geq 0}\|T^n\|^2+1)^{1/2}.$$
We have $$T\cong\begin{pmatrix} V_1 & \ast\\ \mathbb O & T_2\end{pmatrix}.$$ By Lemma \ref{lll02}, 
\begin{equation*}\begin{aligned}\sup_{n\geq 0}\|T^{-n}\| & \leq (\sup_{n\geq 0}\|T_2^{-n}\|)\cdot(2\sup_{n\geq 0}\|T^n\|^2+1)^{1/2} \\ & 
\leq(\sup_{n\geq 0}\|T_1^{-n}\|)\cdot(2\sup_{n\geq 0}\|T^n\|^2+1).\end{aligned}\end{equation*}
\end{proof}

\medskip

{\bf Proof of  Main Theorem.} Denote by $\mathcal H$ the space in which $T$ acts. We have $T=V+(\cdot,v)u$, where $u$, $v\in\mathcal H$ and 
$V\in\mathcal L(\mathcal H)$ is a unitary operator. By Lemma \ref{lllsp} (applied to $T^\ast$), $\sigma(T)\subset\mathbb T$. Thus, $T$ is almost unitary. By Lemma \ref{lllsn} (applied to $T^\ast$), $\mathcal H_s(T^\ast)=\mathcal H$. Therefore, 
$\mathcal H_s(T)=\mathcal H$. By Corollary \ref{cor7.4}, $V$ is a singular unitary operator. 

By Theorem \ref{ttt5}, (\ref{qqq21}) is fulfilled. Without loss of generality, we can suppose that the identity takes place instead of the unitarily 
equivalence in  (\ref{qqq21}).  Denote by $\mathcal H_1$  the space on which $V_1$ acts, and by $X_0$ 
the quasiaffinity which realizes the relation $T\prec U$.  Set $$\mathcal K_1=\operatorname{clos}X_0(\mathcal H_1\oplus L^2(\nu)), \ \   
 \ \ U_1=U|_{\mathcal K_1},$$ 
$$X_1=X_0|_{\mathcal H_1\oplus L^2(\nu)}\in\mathcal L(\mathcal H_1\oplus L^2(\nu), \mathcal K_1).$$
Then $X_1$ is a quasiaffinity, and $$X_1T|_{\mathcal H_1\oplus L^2(\nu)}=U_1X_1.$$ We need the following fact. If a power bounded operator  has the form 
\begin{equation}
\label{bbb}\begin{pmatrix} V_1 & \ast\\ \mathbb O & T_1\end{pmatrix}, 
\end{equation}
where $V_1$ is unitary, then it is similar to $V_1\oplus T_1$, see, 
for example, \cite{b}. Since  $T|_{\mathcal H_1\oplus L^2(\nu)}$ has the form (\ref{bbb}), 
$T|_{\mathcal H_1\oplus L^2(\nu)}\approx V_1\oplus T_1$.  Therefore, there exists an invertible operator $X_2\in\mathcal L(\mathcal H_1\oplus L^2(\nu))$ 
 such that  $$T|_{\mathcal H_1\oplus L^2(\nu)}X_2 = X_2(V_1\oplus T_1).$$ Set $$\mathcal K=\operatorname{clos}X_1X_2 L^2(\nu),  
\ \  U_2=U_1|_{\mathcal K}, \ \ X=X_1X_2|_{L^2(\nu)}\in\mathcal L(L^2(\nu), \mathcal K).$$  Then $X$ is a quasiaffinity, and $XT_1=U_2X$. Recall that $\nu$, $\varphi$, and $\psi$ are defined in  (\ref{qqq21}).
Since $\nu$ is singular (with respect to the Lebesgue measure $m$), and  $\varphi\neq 0$ a.e. with respect to $\nu$, $\varphi$ is cyclic 
for $U_\nu$. By Lemma \ref{lll01}, $\varphi$ is cyclic for $T_1$. Since $T_1\prec U_2$, $U_2$ is cyclic. 

Show that $T_1^\ast$ satisfy the conditions of Proposition \ref{ppp01}  (up to unitarily equivalence). Indeed, 
there exist positive  Borel measures $\mu_\ast$ and $\nu_\ast$ on $\mathbb T$ such that $U_\nu^{-1}\cong U_{\nu_\ast}$ and 
$U_2^{-1}\cong U_{\mu_\ast}$. Note that both $\mu_\ast$ and $\nu_\ast$ are singular with respect to the Lebesgue measure $m$.
  We have  $$T_1^\ast\cong U_{\nu_\ast}+(\cdot,\varphi_\ast)\psi_\ast,$$ 
where $\varphi_\ast$, $\psi_\ast\in L^2(\nu_\ast)$ and $\varphi_\ast\neq 0$, $\psi_\ast\neq 0$ a.e. with respect to $\nu_\ast$, because of 
$\varphi\neq 0$ and  $\psi\neq 0$ a.e. with respect to $\nu$. Set 
$$R_\ast=U_{\nu_\ast}+(\cdot,\varphi_\ast)\psi_\ast,$$ and denote by $Y$ the quasiaffinity 
which realizes the relation $U_{\mu_\ast}\prec R_\ast$. By Corollary \ref{ccc1} applied with $R=U_{\nu_\ast}$, $T=R_\ast$, and $V=U_{\mu_\ast}$,  we obtain that $Y^\ast\varphi_\ast$ is a cyclic vector for $U_{\mu_\ast}^\ast$. Therefore, 
$Y^\ast\varphi_\ast\neq 0$ a.e. with respect to $\mu_\ast$. 

Applying Proposition \ref{ppp01} to $R_\ast$, we obtain  positive  Borel measures $\mu_1$, $\nu_1$ on $\mathbb T$, $\varphi_1\in L^2(\nu_1)$, and a quasiaffinity 
$Z_1\in\mathcal L(L^2(\mu_1),L^2(\nu_1))$ such that $\mu_1(\mathbb T)=1$, $\nu_1(\mathbb T)=1$, $\mu_\ast$ and $\mu_1$ are mutually absolutely continuous, 
$\nu_\ast$ and $\nu_1$ are mutually absolutely continuous, $$Z_1^\ast\overline\chi= \overline\chi, \ \  
Z_1U_{\mu_1}=R_1Z_1, \ \text{  and } \ \ R_\ast\cong R_1, \ \text{ where } \ R_1 = U_{\nu_1}+(\cdot,\overline\chi)\varphi_1.$$
Note that both $\mu_1$ and $\nu_1$ are singular with respect to $m$.

Denote by $\theta$ and $\omega$ the inner functions which are constructed 
by the measures $\mu_1$ and $\nu_1$  as in (\ref{yy0}). Recall that  the unitary operators  $J_{\theta, 1}$ and $J_{\omega, 1}$  are defined  in (\ref{y4}). Set $R=J_{\omega, 1}^{-1}R_1J_{\omega, 1}$ and 
$Z=J_{\omega, 1}^{-1}Z_1J_{\theta, 1}$. Since $$R=U_{(\omega)1}+(\cdot,\overline\chi\omega)u_1 \ \ \text{ for some } \ \ u_1\in\mathcal K_\omega$$ and, by (\ref{y3}), 
$U_{(\omega)1}=S_\omega+(\cdot,\overline\chi\omega)\text{\bf 1}$, we obtain that 
$$R=S_\omega+(\cdot,\overline\chi\omega)u \ \ \text{ for } \ \ u=u_1+\text{\bf 1}\in\mathcal K_\omega.$$ Also, 
$$Z^\ast\overline\chi\omega=\overline\chi\theta, \ \ \ \text{ and } \ \ \ ZU_{(\theta)1}=RZ.$$
By Theorem \ref{ttt3} applied to $R$ and $Z$, there exists $g\in\mathcal K_\omega$ such that 
$Zf=gf$ for $f\in\mathcal K_\theta$. By Theorem \ref{ttt1} applied to $\mu_1$, $R$ and $Z$,  
\begin{equation}
\label{eee}\sup_{n\geq 0}\|R^{-n}\|\leq(\sup_{n\geq 0}\|R^n\|)^5.\end{equation}
 Since $T_1^\ast\cong R_\ast\cong R_1\cong R$, we obtain that (\ref{eee}) is 
fulfilled for $T_1$, too. The conclusion of  Main Theorem follows from (\ref{qqq22}) and (\ref{eee}).

\end{document}